\documentclass[12pt, reqno]{amsart}
\usepackage[english]{babel}
\usepackage{amsmath,amsfonts,amssymb,amsthm}
\usepackage{epsfig}
\usepackage{pgf}
\usepackage{hyperref}
\usepackage{graphicx,graphics}
\usepackage{epstopdf}
\usepackage{color}
\usepackage{caption}
\usepackage{subcaption}
\usepackage{longtable}
\usepackage{underscore}
\usepackage{algorithm}
\usepackage{xfrac}
\usepackage{enumitem}

\usepackage{anyfontsize}

\usepackage[margin=1in]{geometry}
\usepackage[colorinlistoftodos,prependcaption,textsize=tiny]{todonotes}
\definecolor{blue-green}{rgb}{0.0, 0.87, 0.87}
\definecolor{dark-green}{rgb}{0.0, 0.66, 0.00}

\newtheorem{thm}{Theorem}[section]
\newtheorem{innercustomthm}{Theorem}
\newenvironment{customthm}[1]
  {\renewcommand\theinnercustomthm{#1}\innercustomthm}
  {\endinnercustomthm}
\newtheorem{lem}[thm]{Lemma}
\newtheorem{cor}[thm]{Corollary}
\newtheorem{claim}[thm]{Claim}

\newtheorem{problem}[thm]{Problem}
\newtheorem{prop}[thm]{Proposition}
\theoremstyle{definition}
\newtheorem{defin}[thm]{Definition}
\newtheorem{exam}[thm]{Example}
\theoremstyle{remark}

\newcommand { \ib }[1] {\textit{\textbf{#1}}}

\newcommand { \vol }{\mathop{\rm{vol}}\nolimits}
\newcommand { \diag }{\mathop{\rm{diag}}\nolimits}
\newcommand { \trace }{\mathop{\rm{trace}}\nolimits}

\newcommand{\R}{\mathbb{R}}
\newcommand{\Q}{\mathbb{Q}}
\newcommand{\Z}{\mathbb{Z}}
\newcommand{\tr}{\mathrm{tr}}

\begin{document}
\renewcommand{\ib}{\mathbf}
\renewcommand{\proofname}{Proof}
\renewcommand{\phi}{\varphi}
\renewcommand{\epsilon}{\varepsilon}
\newcommand{\conv}{\mathrm{conv}}

\title
[Thinnest lattice sphere covering in dimension six]
{The Vallentin lattice $\mathsf{L}_6^{c}$ gives the thinnest lattice sphere covering in dimension six}

\author[Alexey Garber]{Alexey Garber$^\star$}

\thanks{\thinspace ${\hspace{-.45ex}}^\star$School of Mathematical \& Statistical Sciences,
University of Texas Rio Grande Valley, Brownsville, TX, 78520}

\thanks{\thinspace {\hspace{.45ex}}
Email:
\hskip.06cm
\texttt{alexey.garber@utrgv.edu}}

\date{\today}

\begin{abstract}
We prove that the lattice $\mathsf{L}_6^{c}$ discovered by Vallentin in 2003 and studied by Sch\"urmann and Vallentin in 2006 is the unique (up to similarity) six-dimensional lattice that gives the thinnest lattice sphere covering. 
\end{abstract}

\maketitle

\tableofcontents

\section{Introduction}

The history of questions of optimal placing of spheres goes back to the famous discussion of Newton and Gregory and kissing number in $\R^3$ \cite{PZ04,Zong99} (see also the four-dimensional case \cite{Mus08}) and to Kepler's conjecture, now Hales' theorem \cite{Hal05}, on the optimal packing of spheres in $\R^3$.

Two of the classical problems are related to placing equal spheres centered at points of a lattice in $d$-dimensional Euclidean space, see \cite{Sch09}. Let $v_1,\ldots,v_d$ be a basis in $\R^d$. The lattice $\Lambda$ spanned by the basis is the collection of all integer linear combinations of the basis:
$$\Lambda := \{a_1v_1+\ldots+a_dv_d : a_i\in \Z\}.$$

The lattice sphere packing problem asks to maximize the portion of the space covered by equal nonoverlapping (open) balls centered at the points of the lattice $\Lambda$. The problem is resolved for dimensions $d=2,\ldots,9$ and $d=24$ \cite{Bli35,CK09,DvW25,Gau31,KZ73,KZ77}. In addition to Hales' works on the Kepler conjecture, we also mention works of Viazovska \cite{Via17} and Cohn, Kumar, Miller, Radchenko, and Viazovska \cite{CKMRV17} on optimality of lattice packings using the Korkine-Zolotareff lattice $\mathsf E_8$ and the Leech lattice $\Lambda_{24}$ among all packings of equal spheres in the respective dimensions $d=8$ and $d=24$, respectively.

The second problem is the lattice sphere covering problem that in each dimension asks for the lattice (or lattices) that optimizes covering of the space with equal closed balls centered at the points of the lattice. We start from introducing the main notions and stating the problem rigorously.

For a given lattice $\Lambda$, its {\it covering radius} $R(\Lambda)$ is the smallest radius $R$ so that the (closed) balls of radius $R$ centered at the points of $\Lambda$ cover $\R^d$:
$$R(\Lambda):=\min \{R: \Lambda + R \cdot B^d = \R^d\}$$
where $B^d$ is the closed unit ball in $\R^d$.

The {\it covering density} for the lattice $\Lambda$ is 
$$\Theta(\Lambda):=\frac{R(\Lambda)^d}{\det \Lambda}\kappa_d$$ where $\det \Lambda$ is the fundamental volume of $\Lambda$, i.e. the volume of parallelepiped spanned by any basis of $\Lambda$, and $\kappa_d=\vol (B^d)$. The covering density does not depend on another choice of the basis for the same lattice, or on scaling the lattice, or spanning another lattice with a basis with the same Gram matrix as $v_1,\ldots,v_d$. The latter is equivalent to applying an origin-preserving isometry to the lattice spanned by the basis.

This paper studies the following problem first considered by Kershner in \cite{Ker39}.

\begin{problem}
\label{prob:covering}
For given dimension $d$, find $\Theta_d = \min_\Lambda \Theta(\Lambda)$ where $\Lambda$ runs over all $d$-dimensional lattices and identify all minimizers (up to scaling and isometries).
\end{problem}

The problem is solved in small dimensions while for higher dimensions, candidates to give thinnest lattice coverings are known. The one-dimensional case is trivial. For $d=2$, the problem was solved by Kershner \cite{Ker39} in 1939. Dimension $d=3$ was solved by Bambah \cite{Bam54} in 1954. The four-dimensional case was solved by Delaunay (also spelled Delone)\footnote{Boris Delone (Delaunay), a Russian and Soviet mathematician of French descent. He used the French spelling Delaunay in earlier works, particularly on Delaunay triangulations, and the transliteration of Russian spelling Delone in later works, particularly in works on lattice coverings. We use the French spelling in the paper but the actual spelling in the cited references might be different.} and Ryshkov \cite{DR63}. Finally, the case $d=5$ was solved by Ryshkov and Baranovski in a series of works and was published as a separate monograph \cite{RB78} in 1978. For all dimensions $d=2,3,4,5$, the unique lattice with the thinnest lattice sphere covering is the dual root lattice $\mathsf A_d^*$ where 
$$\mathsf A_d = \{(x_1,\ldots,x_{d+1})\in \Z^{d+1}: x_1+\ldots+x_{d+1}=0\}.$$

For higher dimensions $6\leq d \leq 24$, the candidates for the best sphere coverings are known. For $d=6,7,8$, the corresponding lattices were found by Vallentin \cite{Val03} in 2003 and by Sch\"urmann and Vallentin \cite{SV06} in 2006; the six-dimensional lattice was originally discovered by Vallentin in his PhD dissertation \cite[Sec. 9.6]{Val03} but a more systematic description and conjectured optimality were published in \cite{SV06}. For $9\leq d \leq 21$ except $d=16,18$, the best known lattices were found by Dutour Sikiri\'c, Sch\"urmann, and Vallentin \cite{DSV08} in 2008. The best known lattices in dimensions 16 and 18 are the lattices $\mathsf A_d^*$; their covering densities for all $d$ were computed by Gameckii \cite{Gam62,Gam63} and Bleicher \cite{Ble62}. For $d=24$, the best known lattice is the Leech lattice and its covering density was computed by Conway, Parker, and Sloane \cite{CS88}; the Leech lattice is also known to be locally optimal \cite{SV05} for the lattice sphere covering problem. The lattices for $d=22,23$ were obtained by Smith \cite{Smi88} in 1988 with estimates for their covering densities; see also \cite{DR14} for the precise density of the 23-dimensional lattice. The lattices and the corresponding covering densities are summarized in Table~\ref{table:densities}; we refer to \cite{CS88,SV05,SV06} for details of the lattices mentioned in the table. In higher dimensions $d\geq 24$, Bambah and Sloane \cite{BS82} constructed a series of lattices with better lattice covering densities than the corresponding lattices $\mathsf A_d^*$.

\begin{center}
\begin{table}
{\renewcommand{\arraystretch}{1.2}
\begin{tabular}{r||l|l}
$d$ & Lattice & Covering density\\ \hline
1 & $\mathbb Z$ & 1 \\
2 & $\mathsf A_{2}^{*}$ & 1.209199\\
3 & $\mathsf A_{3}^{*}$ & 1.463503\\
4 & $\mathsf A_{4}^{*}$ & 1.765529\\
5 & $\mathsf A_{5}^{*}$ & 2.124286\\ \hline \hline
6 & $\mathsf{L}_6^{c}$ & 2.464801\\
7 & $\mathsf{L}_7^{c}$ & 2.900024\\
8 & $\mathsf{L}_8^{c}$ & 3.142202\\
9 & $\mathsf{L}_9^{c}$ & 4.268575\\
10 & $\mathsf{L}_{10}^{c}$ & 5.154463\\
11 & $\mathsf{L}_{11}^{c}$ & 5.505591\\
12 & $\mathsf{L}_{12}^{c}$ & 7.465518\\
\end{tabular}
\hskip 1cm
\begin{tabular}{r||l|l}
$d$ & Lattice & Covering density\\ \hline
13 & $\mathsf{L}_{13}^{c}$ & 7.762108\\
14 & $\mathsf{L}_{14}^{c}$ & 8.825210\\
15 & $\mathsf{L}_{15}^{c}$ & 11.004951\\
16 & $\mathsf A_{16}^{*}$ & 15.310927\\
17 & $\mathsf A_{17}^{9}$ & 12.357468\\
18 & $\mathsf A_{18}^{*}$ & 21.840949\\
19 & $\mathsf A_{19}^{10}$ & 21.229200\\
20 & $\mathsf A_{20}^{7}$ & 20.366828\\
21 & $\mathsf A_{21}^{11}$ & 27.773140\\
22 & $\Lambda_{22}^*$ & $\leq 27.8839$\\
23 & $\Lambda_{23}^*$ &  15.3218\\
24 & $\Lambda_{24}$ & 7.903536\\
\end{tabular}}
\caption{Lattices with the least known covering densities for $d \leq 24$. For $d=2,3,4,5$, the lattices are proven to give the best covering density. Our Theorem~\ref{thm:main} proves the six-dimensional case.}
\label{table:densities}
\end{table}
\end{center}

We also mention the currently best known asymptotic upper bound $O(d\ln^\beta d)$ for $\beta\approx 1.858$ for the  lattice sphere covering density which was recently obtained by Gao, Liu, Pikhurko, and Sun \cite{GLPS26} improving on the bound $O(d\ln^\alpha d)$ for $\alpha\approx 2.047$ by Rogers \cite{Rog59} from 1959; see also follow-up by Bukh, Gao, Liu, Pikhurko, and Sun \cite{BGLPS25}.

Our main result in this paper is the proof that the Vallentin lattice $\mathsf{L}_6^c$ is indeed the only lattice with the thinnest lattice sphere covering. This solves Problem~\ref{prob:covering} in the six-dimensional case.

\begin{thm}
\label{thm:main}
The Vallentin lattice $\mathsf{L}_6^{c}$ is the unique (up to similarity) minimum for the lattice sphere covering density in $\R^6$. Particularly, $\Theta_6 = \Theta(\mathsf{L}_6^{c})\approx 2.46480067066519$.
\end{thm}

The rest of the paper is organized as follows. In Section~\ref{sec:pqf} we give necessary definitions to reformulate Problem~\ref{prob:covering} in the language of positive definite quadratic forms and their Delaunay tessellations and triangulations. In Section~\ref{sec:lattice} we give necessary details about the local structure of the Delaunay triangulation of the Vallentin lattice $\mathsf{L}_6^c$ and the auxiliary lattice $\mathsf E_6^*$ that plays a crucial role in our proof. 

Section~\ref{sec:plan} contains an outline of the proof of our main result. In this section we formulate theoretical and computational properties that allow us to complete the proof; details for these statements are provided in the subsequent sections.

In Section~\ref{sec:lemmas} we prove auxiliary lemmas regarding parity classes of vectors, inertia matrices of simplices, and provide details of our computational approach that discards certain local configurations of simplices. Section~\ref{sec:volume} proves that any generic six-dimensional lattice with covering density below $2.4649$ must have a unimodular Delaunay triangulation.

The remaining sections contain details of the computational scripts. In Section~\ref{sec:neighbors}, we restrict local combinatorics of pairs of adjacent simplices in Delaunay triangulations of lattices with small covering densities. In Section~\ref{sec:type-A}, we prove that every simplex of the Delaunay triangulation of a lattice with small covering density must belong to a copy of cluster of 9 simplices observed in the Delaunay triangulation of $\mathsf{L}_6^c$, and in Section~\ref{sec:e6star} we prove that these clusters must repeat the structure of the Delaunay tessellation of $\mathsf E_6^*$. In Section~\ref{sec:exact}, we give details about exact algebraic expressions for lattice $\mathsf L_6^c$ and its Gram matrix $Q_6^c$ and prove its covering density optimality in the corresponding secondary cone.

The paper concludes with several remarks about potential computations in higher dimensions.

\subsection*{Statement of AI use} The author used {\tt Claude} (models {\tt Opus 5.5}, {\tt Fable 5.0}, {\tt Fable 5.1}, and {\tt Fable 5.2}) and {\tt ChatGPT} (models {\tt GPT-5.6 Sol} and {\tt GPT-6.2 Astra}) while working on this paper. Both {\tt Claude} and {\tt ChatGPT} were used to test author's initial proof ideas. {\tt Claude} was used to initially test computational approaches to combinatorial structures used in the proof, and computational scripts for the main Theorem~\ref{thm:main}. {\tt Claude} was used to organize the structure of computational scripts. Both {\tt ChatGPT} and {\tt Claude} were used to obtain and cross-verify the exact expression for the lattice $\mathsf L_6^c$ in Section~\ref{sec:exact}. Both {\tt Claude} and {\tt ChatGPT} were used for literature search, proofreading, and minor editorial suggestions.

Proof ideas for technical Lemmas~\ref{lem:det-inertia},~\ref{lem:trace-inertia},~\ref{lem:int-norm} from Section~\ref{sec:lemmas} and for Lemma~\ref{lem:vol1} were initially suggested by {\tt Claude} and adapted and simplified by the author.

All proofs were written and verified by the author. All computational scripts were written from scratch and verified by the author using {\tt Claude}'s help. The {\tt SageMath} script to verify the exact form $Q_6^c$ was written by the author. The author is responsible for all proofs and computational details of this paper.

\subsection*{Acknowledgements} Author thanks Frank Vallentin, Sven Ahrend, and Achill Sch\"urmann for helpful comments on an early version of this paper.

After the initial version of this work was made publicly available, the author learned about an alternative and independent approach by Ahrend and Sch\"urmann to proving Theorem~\ref{thm:main}; their manuscript~\cite{AS26} is currently in preparation.

Author's research was partially supported by Simons Foundation.

\section{Covering radius through positive definite quadratic forms}
\label{sec:pqf}

\subsection{Voronoi polytopes and Delaunay tessellation}

Let $\Lambda$ be a $d$-dimensional lattice. We first note that the {\it Dirichlet-Voronoi polytope} of $\Lambda$ consisting of all points of $\R^d$ that are closer to the origin than to any other point of $\Lambda$ is helpful to find the covering radius of $\Lambda$. Translated copies of the Dirichlet-Voronoi polytope tile $\R^d$ and consequently, the covering radius $R(\Lambda)$ is the maximal distance from vertices of the Dirichlet-Voronoi polytope to the origin. A related approach is to consider the tiling dual to the tiling by Dirichlet-Voronoi polytopes.

\begin{defin}
\label{def:delaunay}
For a given $d$-dimensional lattice $\Lambda$, the {\it Delaunay tessellation} of $\Lambda$, denoted $\mathcal D(\Lambda)$ is the tiling of $\R^d$ with the vertex set $\Lambda$ that consists of all ($d$-dimensional) polytopes with vertices in $\Lambda$ that are inscribed in spheres that do not contain points of $\Lambda$ inside.

If $\mathcal D(\Lambda)$ consists of simplices only, then we also call it {\it Delaunay triangulation} of $\Lambda$.
\end{defin}

It is known that $\mathcal D(\Lambda)$ is a face-to-face tiling of $\R^d$. The tessellation can be obtained using the lifting to paraboloid $f:\R^d \longrightarrow \R^{d+1}$ where $f(x)=(x,||x||^2)$ and $||x||$ is the standard Euclidean norm. Then the Delaunay tessellation  of $\Lambda$ is the projection of the boundary of the convex hull of $f(\Lambda)$ back to $\R^d$; note that the convex hull is unbounded in the positive direction of $x_{d+1}$ .

The covering radius $R(\Lambda)$ in this setting is the largest radius among circumradii of the polytopes in the Delaunay tessellation $\mathcal D(\Lambda)$.

\subsection{Positive definite quadratic forms}

We will work in a dual setting where every lattice can be identified with the integer lattice and different Delaunay triangulations are obtained by considering different Euclidean metrics in $\R^d$. Let $v_1,\ldots,v_d$ be a basis of $\Lambda$ and $e_1,\ldots,e_d$ be the standard orthonormal basis in $\R^d$. The linear transformation that maps $v_i \mapsto e_i$ maps the lattice $\Lambda$ to the standard integer lattice $\Z^d$. We can recover the image of the Delaunay tessellation of $\Lambda$ by considering the squared norm $Q[x]:=x^tQx$ associated with the Gram matrix $Q$ of the basis $v_1,\ldots,v_d$ and the corresponding metric $d_Q(x,y)=\sqrt{Q[x-y]}$ on $\R^d$ induced by that norm.

Thus, instead of parameterizing the lattices with their bases, we will consider another parameterization of Problem~\ref{prob:covering} using positive definite quadratic forms $Q$ in the convex cone $\mathcal S_d^{>0}$ of all positive definite quadratic forms.

\begin{defin}
Let $Q$ be a positive definite quadratic form, i.e. a positive definite $d\times d$ symmetric matrix $Q$. The {\it Delaunay tessellation} (of $\Z^d$) associated with $Q$ is the collection of all polytopes with vertices in $\Z^d$ that are inscribed in empty ellipsoids of the form $(x-c)^tQ(x-c)=a^2$.

If a Delaunay tessellation consists of simplices only, then we call it a {\it Delaunay triangulation}.
\end{defin}

For each ellipsoid of the form $(x-c)^tQ(x-c)=a^2$, we call $c$ its {\it center}, and (positive) $a$ its {\it radius}.

\begin{defin}
For a given positive definite quadratic form $Q$, its {\it covering radius} $R(Q)$ is the minimum among all positive $a$ such that the family of (solid) ellipsoids $$\{(x-c)^tQ(x-c)\leq a^2: c\in \Z^d\}$$
covers $\R^d$. 
\end{defin}

Similarly to the standard setting, the covering radius $R(Q)$ is the maximum among radii of circumscribed ellipsoids defined by $Q$ for the polytopes in the Delaunay tessellation $\mathcal D(Q)$.

The covering density of the quadratic form $Q$ is then defined as 
$$\Theta(Q)=\frac{R(Q)^d}{\sqrt{\det Q}}\kappa_d$$
and Problem~\ref{prob:covering} is equivalent to finding the minimum  $\Theta_d=\min\{\Theta (Q):Q\in \mathcal S_d^{>0}\}$ and finding all positive definite quadratic forms that achieve this minimum.

\subsection{Secondary cones}

One additional thing we will need in this setting is the notion of $L$-type domain or secondary cone introduced by Voronoi \cite{Vor08} and developed by Delaunay \cite{Del38}.

\begin{defin}
Let $\mathcal D$ be the Delaunay triangulation of $\Z^d$ associated with a positive definite quadratic form. The {\it secondary cone} of $\mathcal D$, also known as the {\it $L$-type domain} of $\mathcal D$ is the collection of all positive definite quadratic forms $Q$ such that $\mathcal D(Q)=\mathcal D$. That is,
$$SC(\mathcal D)=\{Q\in \mathcal S_d^{>0}: \mathcal D(Q)=\mathcal D\}.$$
\end{defin}

It is not complicated to see that each secondary cone is a convex subcone in the cone $\mathcal S_d^{>0}$ of all positive definite quadratic forms. We refer to \cite{DGSW16} for the following properties of secondary cones.

\begin{lem}
\label{lem:closure}
Let $SC(\mathcal D)$ be the secondary cone of a Delaunay triangulation $\mathcal D$.
\begin{enumerate}
\item The (rational) closure of $SC(\mathcal D)$ is a convex polyhedral cone in the cone of positive semidefinite quadratic forms.
\item A positive definite quadratic form $Q'$ is in the closure of $SC(\mathcal D)$ if and only if $\mathcal D$ is a subdivision of $\mathcal D(Q')$.
\end{enumerate}
\end{lem}

The study of equivalence classes of secondary cones is at the heart of the Voronoi second reduction theory for lattices (or positive definite quadratic form) and we refer to \cite{DGSW16} for more details. Another problem closely related to secondary cones is the Voronoi conjecture \cite{Gar25} that claims that every convex polytope that tiles $\R^d$ with translations can be obtained as the Dirichlet-Voronoi polytope of a lattice with respect to some positive definite quadratic form. The study of such polytopes, called {\it parallelohedra}, plays a prominent role in mathematical crystallography and goes to works of Fedorov \cite{Fed85} who classified three-dimensional parallelohedra and works of Minkowski \cite{Min97} that was one of the inspirations for geometry of numbers.

{\bf Notations.} We use $e_1,\ldots,e_d$ for the vectors of the standard orthonormal basis in $\R^d$.

Once we fix the quadratic form $Q$, the terms as radius, circumradius, circumcenter will refer to the metric defined by $Q$. In particular, we will say that the radius of ellipsoid $(x-c)^tQ(x-c)=a^2$ is $a$ (provided it is positive). The circumcenter of the (full-dimensional) simplex $L$ is the center $c$ of the ellipsoid $(x-c)^tQ(x-c)=a^2$ that contains all vertices of $L$ and its circumradius is the radius of this ellipsoid. Similarly, lengths will be referred with respect to $Q$. So squared length, or squared norm, of segment, or vector, $v$ is $v^tQv$.

We also use normalized volume for simplices involved. If $S=\conv\{v_0,\ldots,v_d\}\subset \R^d$, then the {\it normalized volume} of $S$ is the absolute value of the determinant $\det(v_1-v_0,\ldots,v_d-v_0)$. The reason for using the normalized volume is to use integer volumes for simplices with integer coordinates rather than volumes with denominator $d!$.

\section{The Vallentin lattice \texorpdfstring{$\mathsf{L}_6^{c}$}{L6c}}
\label{sec:lattice}

In this section we describe lattice $\mathsf L_6^c$, the associated positive definite quadratic form $Q_6^c$, and the corresponding Delaunay triangulation. We start from the lattice $\mathsf E_6^*$ and its definition through the Korkine-Zolotareff lattice $\mathsf E_8$. We refer to \cite{CS91} and \cite{SV06} for details of the constructions. We also note that we make certain choices to obtain the setting of Sch\"urmann and Vallentin \cite{SV06} but our further computational steps will use a slightly different setup which is equivalent to the constructions in this section via a $\Z^6$-preserving affine transformation.

\subsection{Lattice \texorpdfstring{$\mathsf E_6^*$}{E6*}}

The lattice $\mathsf E_8$ is defined as 
$$\mathsf E_8 = \left\{(x_1,\ldots,x_8) \in \Z^8 \cup \left(\left(\frac12,\frac12,\frac12,\frac12,\frac12,\frac12,\frac12,\frac12\right)^t+\Z^8\right): \sum x_i \text{ is even}\right\}.$$
The lattice $\mathsf E_8$ has $240$ minimal vectors of squared length $2$.

The lattice $\mathsf E_6$ is a slice of $\mathsf E_8$ orthogonally to a hexagonal sublattice generated by two shortest vectors of $\mathsf E_8$. For example, we can take $\mathsf E_6$ to consist of vectors of $\mathsf E_8$ in the six-dimensional subspace $x_6=x_7=x_8$ orthogonal to minimal vectors $(0,0,0,0,0,0,1,-1)^t$ and $(0,0,0,0,0,1,0,-1)^t$ of $\mathsf E_8$. The following six vectors form a basis for the lattice $\mathsf E_6$:

$$\begin{array}{l}
u_1 = (1, 1, 0, 0, 0, 0, 0, 0)^t\\
u_2 = (1, -1, 0, 0, 0, 0, 0, 0)^t\\
u_3 = (-1, 0, -1, 0, 0, 0, 0, 0)^t\\
u_4 = (-1, 0, 0, -1, 0, 0, 0, 0)^t\\
u_5 = (0, 1, 0, 0, -1, 0, 0, 0)^t\\
u_6 = \left(\frac12, -\frac{1}2, \frac12, \frac12, \frac12, -\frac{1}2, -\frac{1}2, -\frac{1}2\right)^t.
\end{array}$$

The Gram matrix of this basis is $$
Q_{\mathsf E_6}=
\left(\begin{array}{cccccc}
2 & 0 &-1 &-1&  1 & 0\\
0 & 2 &-1& -1& -1&  1\\
-1& -1&  2&  1&  0& -1\\
-1& -1&  1&  2&  0& -1\\
1& -1&  0&  0&  2& -1\\
0&  1& -1& -1& -1&  2
\end{array}\right)
$$
and therefore the Gram matrix of the dual basis (after appropriate scaling), and the positive definite quadratic form associated with $\mathsf E_6^*$, is 
$$Q_{\mathsf E_6^*}=\left(
\begin{array}{cccccc}
4 & 1 & 2 & 2 & -1 & 1 \\    
1 & 4 & 2 & 2 & 2  & 1 \\
2 & 2 & 4 & 1 & 1  & 2 \\
2 & 2 & 1 & 4 & 1  & 2 \\
-1& 2 & 1 & 1 & 4  & 2 \\
1 & 1 & 2 & 2 & 2  & 4
\end{array}
\right).$$

The Delaunay tessellation of the quadratic form $Q_{\mathsf E_6^*}$, and of the lattice $\mathsf E_6^*$ consists of 40 classes of free sums of three triangles. More precisely, let $\mathcal P$ be a convex hull of three regular triangles with a common center $O$ located in three pairwise orthogonal two-dimensional subspaces spanning the six-dimensional space. There are 40 polytopes $P_1,\ldots, P_{40}$ affinely equivalent to $\mathcal P$ and the Delaunay tessellation of $\mathsf E_6^*$ consists of translated copies of the polytopes $P_i$ for $i=1,\ldots,40$ and their centrally symmetric copies. We refer to \cite{CS91} for more details including description that uses the corresponding symmetry group of the lattice.

For our purposes, we will consider potential triangulations of two free sums adjacent by a facet. Let $P_1$ and $P_2$ be the polytopes of the Delaunay tessellation of $\Z^6$ with respect to $Q_{\mathsf E_6^*}$ with circumcenters $\frac13(1,1,0,0,0,1)^t$ and $\frac13(1,0,1,1,1,-1)^t$; the circumradii are $\sqrt2$ for both polytopes. More precisely,
\begin{equation}
\label{eq:P1}
P_1 = \conv \left(\begin{array}{ccc}
(0,0,0,0,0,0)^t & (1,0,0,0,1,0)^t & (0,1,0,0,-1,1)^t \\
(1,0,0,0,0,0)^t & (0,1,0,0,0,0)^t & (0,0,0,0,0,1)^t \\
(0,0,1,0,0,0)^t & (0,0,0,1,0,0)^t & (1,1,-1,-1,0,1)^t \\
\end{array} \right)
\end{equation}
and 
\begin{equation}
\label{eq:P2}
P_2 = \conv \left(\begin{array}{ccc}
(0,0,0,0,0,0)^t & (1,0,0,0,1,0)^t & (0,0,1,1,0,-1)^t \\
(1,0,0,0,0,0)^t & (0,1,0,0,0,0)^t & (0,-1,1,1,1,-1)^t \\
(0,0,1,0,0,0)^t & (0,0,0,1,0,0)^t & (1,0,0,0,1,-1)^t \\
\end{array} \right).
\end{equation}

\subsection{Lattice \texorpdfstring{$\mathsf L_6^c$}{L6c}}

The quadratic form $Q_6^c$ associated with the lattice $\mathsf L_6^c$ has Delaunay triangulation that subdivides the Delaunay tessellation of $\mathsf E_6^*$; consequently, the quadratic form $Q_{\mathsf E_6^*}$ is in the closure of the secondary cone that contains $Q_6^c$ by Lemma~\ref{lem:closure}. Our next goal is to describe potential neighbors in all such triangulations including the triangulation corresponding to $\mathsf L_6^c$. It is enough to work only with the polytopes $P_1$ and $P_2$ as any two pairs of polytopes in the Delaunay tessellation of $\mathsf E_6^*$ with common facet can be matched by the symmetry group of $\mathsf E_6^*$.

\begin{defin}
Let $\sigma_1=\conv\{v_0,\ldots,v_6\}$ and $\sigma_2=\conv\{v_1,\ldots,v_7\}$ be two unimodular simplices with integer vertices in $\Z^6$ with a common facet $F=\conv\{v_1\ldots,v_6\}$. Suppose $v_0$ and $v_7$ are on different sides of $F$.

Let $x$ be the coordinates of the vector $v_7-v_0$ in the basis $v_1-v_0,\ldots,v_6-v_0$ sorted in the non-increasing order. We call $x$ the {\it type} of the facet $F$.
\end{defin}

We note that the type of a facet is defined by the pair of simplices sharing the facet not only by the facet itself. Moreover, switching the order of simplices does not change the type of $F$ as it can be reconstructed from the barycentric coordinates of the point $\frac12(v_0+v_7)$ with respect to the vertices of $F$. We also note that each type is just a non-increasing sequence of six integers adding up to $2$.

We will be most interested in the following three types of facets.
\begin{center}
\begin{tabular}{ll}
Type A:& $(1,1,1,0,0,-1)$\\
Type B:& $(1,1,0,0,0,0)$\\
Type C:& $(1,1,1,1,-1,-1)$
\end{tabular}
\end{center}

The following lemma shows which local properties we will be pursuing in our search of local collection of simplices within potentially best lattice for sphere covering in $\R^6$.

\begin{lem}
\label{lem:typesABC}
Let $\sigma_1$ and $\sigma_2$ be two simplices of a refinement of the Delaunay tessellation corresponding to $Q_{\mathsf E_6^*}$ with a common facet $F$.
\begin{enumerate}
    \item If $\sigma_1$ and $\sigma_2$ are in one free sum of three triangles from the Delaunay tessellation of $Q_{\mathsf E_6^*}$, then $F$ is of Type A.
    \item If $\sigma_1$ and $\sigma_2$ are in two different free sums of three triangles from the Delaunay tessellation of $Q_{\mathsf E_6^*}$, then $F$ is of Type B or of Type C.
\end{enumerate}
\end{lem}

\begin{proof}
Suppose $\sigma_1$ and $\sigma_2$ are in one free sum of three triangles; we may assume this is the sum $P_1$, see \eqref{eq:P1}. There are three triangulations of $P_1$ with 9 simplices each. For these triangulations, we pick one of the triangles (three points in one row of \eqref{eq:P1}), and nine combinations of two points from each of the other two triangles. For a given simplex, an adjacent simplex can be obtained by swapping a point from a non-picked triangle. For all three triangulations, any two pairs with common facets are affinely equivalent and the type can be verified for any pair. 

For example, if
$$
\sigma_1 = \conv \left(\begin{array}{c}
(0,0,0,0,0,0)^t \\ (1,0,0,0,1,0)^t \\ (0,1,0,0,-1,1)^t \\
(1,0,0,0,0,0)^t \\ (0,1,0,0,0,0)^t \\ 
(0,0,1,0,0,0)^t \\ (0,0,0,1,0,0)^t \\ 
\end{array} \right)
\quad
\text{and}
\quad
\sigma_2 = \conv \left(\begin{array}{c}
(0,0,0,0,0,0)^t \\ (1,0,0,0,1,0)^t \\ (0,1,0,0,-1,1)^t \\
(0,0,0,0,0,1)^t \\ (0,1,0,0,0,0)^t \\ 
(0,0,1,0,0,0)^t \\ (0,0,0,1,0,0)^t \\ 
\end{array} \right)
$$
with vertex $(1,0,0,0,0,0)^t$ of $\sigma_1$ replaced by vertex $(0,0,0,0,0,1)^t$ of $\sigma_2$, then 
{\tiny
$$
\left(
\begin{array}{c}
0\\0\\0\\0\\0\\1
\end{array}
\right)-
\left(
\begin{array}{c}
1\\0\\0\\0\\0\\0
\end{array}
\right)
=
\left(
\left(
\begin{array}{c}
0\\0\\0\\0\\0\\0
\end{array}
\right)-
\left(
\begin{array}{c}
1\\0\\0\\0\\0\\0
\end{array}
\right)
\right) 
+
\left(
\left(
\begin{array}{c}
1\\0\\0\\0\\1\\0
\end{array}
\right)-
\left(
\begin{array}{c}
1\\0\\0\\0\\0\\0
\end{array}
\right)
\right)
+
\left(
\left(
\begin{array}{c}
0\\1\\0\\0\\-1\\1
\end{array}
\right)-
\left(
\begin{array}{c}
1\\0\\0\\0\\0\\0
\end{array}
\right)
\right)
-
\left(
\left(
\begin{array}{c}
0\\1\\0\\0\\0\\0
\end{array}
\right)-
\left(
\begin{array}{c}
1\\0\\0\\0\\0\\0
\end{array}
\right)
\right)
$$
}
and the common facet of these two simplices is of Type A.

If $\sigma_1$ and $\sigma_2$ are in different sums of three triangles, then we can assume that $\sigma_1$ is in $P_1$ from \eqref{eq:P1} and $\sigma_2$ is in $P_2$ from \eqref{eq:P2} while the common facet is $P_1\cap P_2$. There are three options for the remaining vertex of $\sigma_1$ and three options for the remaining vertex of $\sigma_2$ and we can verify the types for all 9 pairs. Alternatively, we can use symmetries of $P_1$ fixing the common facets, and verify only the pairs with one fixed $\sigma_1$. In all cases, the computations can be checked in a similar way as in the first case.
\end{proof}

Now we are going to describe the quadratic form $Q_6^c$ corresponding to the Vallentin lattice $\mathsf L_6^c$. The lattice is known only approximately and we describe the precise quadratic form using roots of three polynomials of degrees 19 in Section~\ref{sec:exact}.

Let $R_2$, $R_3$, and $R_4$ be the three positive semidefinite quadratic forms given by
{\small$$
R_2=\left(\begin{array}{cccccc}
0& 0& 0& 0& 0& 0\\
0& 5& 2& 2& 3& 1\\
0& 2& 4& 0& 2& 2\\
0& 2& 0& 4& 2& 2\\
0& 3& 2& 2& 5& 3\\
0& 1& 2& 2& 3& 5
\end{array}\right),
\quad
R_3=\left(\begin{array}{cccccc}
6& 4& 4& 4& 0& 2\\
4& 11& 6& 6& 5& 3\\
4& 6& 8& 3& 3& 4\\
4& 6& 3& 8& 3& 4\\
0& 5& 3& 3& 7& 4\\
2& 3& 4& 4& 4& 7
\end{array}\right),
\quad
R_4=\left(\begin{array}{cccccc}
3& 2& 2& 2& 0& 1\\
2& 3& 2& 2& 1& 1\\
2& 2& 4& 1& 1& 2\\
2& 2& 1& 4& 1& 2\\
0& 1& 1& 1& 3& 2\\
1& 1& 2& 2& 2& 4
\end{array}\right).
$$}
Then 
\begin{equation}
\label{eq:approx}
Q_6^c = Q_{\mathsf E_6^*} + x R_2 + y R_3 + z R_4
\end{equation}
where $x = 0.15266328480099, y =
0.32884740614948, z = 0.13827491241153$. We note that the values $x,y,z$ are approximate and Sch\"urmann and Vallentin found these values by solving the lattice covering problem over the secondary cone of $Q_6^c$. The quadratic forms $Q_{\mathsf E_6^*}$, $R_2$, $R_3$, and $R_4$ are on the boundary of this cone and the linear combination $Q_6^c$ is strictly inside the cone. We refer to \cite{SV06} for more details.

The most important property of the lattice $\mathsf L_6^c$ and the corresponding quadratic form $Q_6^c$ that we are going to use is that it solves the lattice sphere covering in a neighborhood of $\mathsf E_6^*$ within the cone $\mathcal S_6^{>0}$.

\begin{thm}[Sch\"urmann, Vallentin, {\cite[Thm. 11.6]{SV06}}]
\label{thm:SV}
Among all positive definite quadratic forms with Delaunay tessellations subdividing the Delaunay tessellation of $Q_{\mathsf E_6^*}$, the quadratic form $Q_6^c$ gives the least dense lattice sphere covering.
\end{thm}

We note that Sch\"urmann and Vallentin actually proved a stronger statement. We formulate it here as proposition and it follows from the proof of \cite[Thm. 11.6]{SV06}.

\begin{prop}[Sch\"urmann, Vallentin, \cite{SV06}]
\label{prop:SV-unique}
Suppose $Q$ is a positive definite quadratic form such that the Delaunay triangulation $\mathcal D(Q)$ subdivides the Delaunay tessellation of $\mathsf E_6^*$. If $\Theta(Q)\leq 2.464802$, then $Q$ belongs to a secondary cone $GL_6(\Z)$-equivalent to the secondary cone of $Q_6^c$ given by \eqref{eq:approx}.
\end{prop}

\section{Proof of the main theorem}
\label{sec:plan}

In this section we outline the proof of our main theorem. Our goal is to prove that $\mathsf{L}_6^{c}$ is the unique lattice (up to scaling) minimizing the density of lattice sphere covering in $\R^6$. Equivalently, that the corresponding quadratic form $Q_6^{c}$ is the only positive definite quadratic form (up to scaling and $GL_6(\Z)$-equivalence) minimizing the covering density $\Theta(Q)$. Suppose $Q$ is another form with $\Theta(Q)\leq \Theta(Q_6^{c})\leq 2.464801$.

For now, let's assume that the Delaunay tessellation corresponding to the quadratic form $Q$ is a triangulation so it is inside the corresponding full-dimensional secondary cone; we'll treat the boundary case in the proof of our main theorem below. 

We will also assume that the covering radius $R(Q)$ of $Q$ is $1$ as scaling does not change covering density. And since $\Theta(Q)=\frac{R(Q)^6}{\sqrt{\det Q}}\kappa_6$ with $\kappa_6=\frac{\pi^3}{6}$ being the volume of the six-dimensional unit ball, we get 
$$\det(Q)=\frac{\kappa_6^2}{\Theta(Q)^2}>4.395397.$$

We will find upper bounds on $\det(Q)$ from potential combinatorial structures within its Delaunay triangulation. If such a bound is less than $4.395397$, then it implies that no Delaunay triangulation of $\Z^6$ that contains such a combinatorial structure defines a secondary cone that contains a quadratic form with smaller covering density than $Q_6^{c}$.

Below we provide a proof of our main theorem that is based on volume arguments from Section~\ref{sec:volume} and computational statements that are established in Sections~\ref{sec:neighbors},~\ref{sec:type-A}, and~\ref{sec:e6star}. The statements make sure that any generic $Q$ with $\det Q > 4.395397$ and covering radius at most~$1$ must have local structure of its Delaunay triangulation resembling the one of the Delaunay triangulation of $Q_6^{c}$ or any subdivision of the Delaunay tessellation of $Q_{\mathsf E_6^*}$ described in Lemma~\ref{lem:typesABC}. After that, we use the optimality of $Q_6^c$ for covering density among quadratic forms subdividing Delaunay tessellation of $\mathsf E_6^*$ shown by Sch\"urmann and Vallentin and Lemma~\ref{lem:exact} that shows that $Q_6^c$ is unique optimizer in its secondary cone.

\begin{customthm}{\ref{thm:main}}
The Vallentin lattice $\mathsf{L}_6^{c}$ is the unique (up to similarity) minimum for the lattice sphere covering density in $\R^6$. Particularly, $\Theta_6 = \Theta(\mathsf{L}_6^{c})\approx 2.46480067066519$.
\end{customthm}

\begin{proof}
We will prove the theorem in the language of positive definite quadratic forms. It is enough to show that $Q_6^c$ is the only quadratic form minimizing the covering density up to scaling and $GL_6(\Z)$ transformations. Suppose there is a six-dimensional positive definite quadratic form $Q$ with $\Theta (Q) \leq \Theta(\mathsf L_6^c)$. 

First, let's assume that the Delaunay tessellation corresponding to the form $Q$ is a triangulation. According to Theorem~\ref{thm:SV} and Proposition~\ref{prop:SV-unique}, it is enough to show that Delaunay tessellation of $Q$ subdivides the Delaunay tessellation of $Q_{\mathsf E_6^*}$ (or its $GL_6(\Z)$-equivalent) and that the quadratic form $Q_6^c$ is the unique minimizer of the covering density within its secondary cone. 

For the first part, from the discussion right before the theorem, we may assume that all simplices of $\mathcal D(Q)$ have circumradii at most 1, and get $\det Q > 4.395397$ and it is still enough to show that the Delaunay triangulation of $Q$ subdivides the Delaunay tessellation of $Q_{\mathsf E_6^*}$.

We establish this property with the following steps.

In Lemma~\ref{lem:vol1}, we show that the Delaunay triangulation of $Q$ is unimodular. That is, every simplex of the Delaunay triangulation of $Q$ with vertices in $\Z^6$ has normalized volume 1.

In Lemmas~\ref{lem:5types} and~\ref{lem:3types}, we use the unimodularity to show that all facets of the Delaunay triangulation of $Q$ are of Types A, B, or C. This is done by first discarding all but five types that include Types D and E below in addition to the three types in Lemma~\ref{lem:typesABC}.
\begin{center}
\begin{tabular}{ll}
Type D:& $(2,1,1,0,-1,-1)$\\
Type E:& $(1,1,1,1,0,-2)$
\end{tabular}
\end{center}
After that, we run a Depth-First Search in Lemma~\ref{lem:3types} over potential simplices we can further attach to a corresponding pair with common facet of Type D or Type E. For both initial pairs we always reach a collection $\mathcal L$ of simplices that guarantee that $Q$ has determinant less than $4.395397$ if all simplices in $\mathcal L$ have circumradii at most 1.

At this point, we have established that $\mathcal D(Q)$ is unimodular and has only three types of facets. We assume that the standard simplex $\sigma_0$ with vertices at the origin and at all six basis vectors is a simplex of $\mathcal D(Q)$. In Claims~\ref{claim1} and~\ref{claim2}, and in Corollary~\ref{cor:free-sum}, we show that if one facet of $\sigma_0$ in $\mathcal D(Q)$ is of Type A, then we can find nine simplices (including $\sigma_0$) giving a triangulation of a free sum of three triangles. In this subset of nine simplices, each facet common to two simplices is of Type A, cf. Lemma~\ref{lem:typesABC}.

In Claim~\ref{claim3}, we show that once we found this collection of nine simplices, the remaining three neighbors of $\sigma_0$ are of Types B or C (so not of Type A). This shows that if we consider connected components of the graph with the vertex set given by simplices of $\mathcal D(Q)$ and edges given by Type A facets, then simplices of each connected component form a triangulation of a free sum of three triangles.

In Claim~\ref{claim4}, we show that every simplex in $\mathcal D(Q)$ must have at least one Type A facet. The previous property implies that all simplices of $\mathcal D(Q)$ split into subsets of nine simplices forming free sums of three triangles each.

The remaining step is to show that these free sums of three triangles form the Delaunay tessellation of $Q_{\mathsf E_6^*}$ or an equivalent quadratic form. This is done in Claim~\ref{claim-adj-sums}. In this claim we show that one cluster of nine simplices in the Delaunay triangulation $\mathcal D(Q)$ forming a free sum of three triangles, defines the neighboring free sums and the neighbors are the same for all three potential triangulations of the starting free sum. While this property does not give a unique triangulation of $\R^6$ with the vertex set $\Z^6$, it shows that $\mathcal D(Q)$ subdivides the Delaunay tessellation of $Q_{\mathsf E_6^*}$ because the corresponding clusters must be the same as the clusters for $Q_6^c$ (which satisfies all the properties we use including the determinant bound).

After that we can apply Theorem~\ref{thm:SV} and Proposition~\ref{prop:SV-unique} because the Delaunay triangulation $\mathcal D(Q)$ subdivides Delaunay tessellation of $Q_{\mathsf E_6^*}$. It follows that $Q$ is in the secondary cone of $Q_6^c$ and since $\Theta(Q)\leq \Theta(Q_6^c)$, Lemma~\ref{lem:exact} implies that $Q$ is proportional to $Q_6^c$. 

Now we consider the case when the Delaunay tessellation of $Q$ is not a triangulation. Since $\Theta(Q)\leq \Theta(Q_6^c)$ and Sch\"urmann and Vallentin bounded the covering density of $Q_6^c$ from above by $2.464801$, we can perturb into a form $Q_{gen}$ with Delaunay triangulation and with covering density $\Theta(Q_{gen})<2.464802$, cf. Proposition~\ref{prop:SV-unique}. Then our approach to triangulations from this proof implies that $Q_{gen}$ is in the secondary cone of $Q_6^c$ or an equivalent cone. Since we can make the perturbation arbitrarily small, the form $Q$ is on the boundary of the secondary cone of $Q_6^c$ or an equivalent cone. Now we can apply Lemma~\ref{lem:exact} and use that $Q_6^c$ is the only covering density optimizer in the closure of the cone so $Q$ is a scaling of $Q_6^c$.

For the exact covering density, we refer to computations in Section~\ref{sec:exact}.
\end{proof}

\section{Auxiliary lemmas}
\label{sec:lemmas}

In this section we introduce several auxiliary notions and establish several properties that will be used in the proofs of lemmas we referred to in the proof of Theorem~\ref{thm:main}. We also give details on how convex optimization is used to get certified bounds for determinant of our potential positive definite quadratic forms.

We restrict ourselves to the six-dimensional case. So all vectors are with six coordinates and quadratic forms are represented by $6\times 6$ matrices unless stated otherwise.

\subsection{Parity classes}

We start by highlighting restrictions on which sets of vectors may serve as sets of edges of simplices in Delaunay triangulations of positive definite quadratic forms. As a simple observation, if $v\in \Z^6$ has only even entries, then it cannot belong to any Delaunay triangulation or Delaunay tessellation because any ellipsoid that contains its endpoints must contain its midpoint strictly inside.

\begin{defin}
Let $u,v \in \Z^6$ be two integer vectors. We say that $u$ and $v$ belong one {\it parity class} if $u-v$ is a vector with even coordinates. Equivalently, if $u-v\in 2\Z^6$.
\end{defin}

Overall, there are $63=2^6-1$ non-trivial classes represented by non-zero vertices of the standard unit cube $[0,1]^6$. We need the following well-known property of the representatives of parity classes in Delaunay triangulations; see \cite{Vor08}.

\begin{lem}
\label{lem:parity}
Let $Q$ be a positive definite quadratic form with Delaunay triangulation $\mathcal D(Q)$. Then
\begin{enumerate}
\item Each non-zero parity class is represented by an edge of $\mathcal D(Q)$.
\item If two edges of $\mathcal D(Q)$ represent two vectors $u,v \in \Z^6$ from the same parity class, then $u=\pm v$.
\end{enumerate}
\end{lem}

\begin{proof}
A short idea for the proof is the following. For each parity class, consider the shortest vector within this parity class with respect to $Q$. This vector will be the unique vector of the parity class represented in $\mathcal D(Q)$.
\end{proof}

Our main idea to use parity classes will be as follows. Suppose we have a finite collection of simplices $\{\sigma_1,\ldots,\sigma_k\}$ that can potentially belong to a Delaunay triangulation $\mathcal D(Q)$ of some positive definite quadratic form $Q$. If we will try to add one additional simplex $\tau$ to that configuration, then the edges of $\tau$ that represent parity classes already present among edges of $\{\sigma_1,\ldots,\sigma_k\}$ must be parallel to these present edges. We will see how this restricts the potential new simplices that can be added to this collection.

\subsection{Inertia matrix}

Our second tool is the inertia matrix of a simplex with respect to a given positive definite quadratic form. Moments of inertia were used by Sch\"urmann and Vallentin \cite{SV06} to obtain local lower bounds for covering density. We will use the eigenvalues of inertia matrix to restrict types of facets in the Delaunay triangulation of $Q$ based on a lower bound for $\det Q$.

Let $Q$ be a positive definite quadratic form and let $A$ be any $6\times 6$ matrix such that $Q=A^tA$ (for example, we can take $A=Q^{\frac12}$). Let $L=\conv\{v_0,\ldots,v_6\} \subset \R^6$ be a six-dimensional simplex with centroid $m=\frac17(v_0+\ldots+v_6)$; we use column vectors in this section.

\begin{defin}
\label{def:inertia}
We will call the $6\times 6$ matrix 
$$M_L=M_L(A)=\sum_{k=0}^6 A(v_k-m)(v_k-m)^t A^t$$
the {\it inertia matrix} for the simplex $L$ with respect to $Q$.
\end{defin}

We note that if $Q=B^tB$ for another $6\times 6$ matrix $B$, then $U=AB^{-1}$ is an orthogonal matrix and inertia matrix $M_L(B)=U^tM_L(A)U$. So the inertia matrix is defined up to conjugation with an orthogonal matrix. In particular, eigenvalues of $M_L$ do not depend on the choice of $A$.

We also add that $M_L$ is a positive semidefinite matrix because it is a sum of positive semidefinite matrices so the eigenvalues of $M_L$ are nonnegative real numbers. Our goal is to use spectral properties of $M_L$ to obtain bounds on lengths of integer vectors with respect to $Q$. We start from the determinant.

\begin{lem}
\label{lem:det-inertia}
Let $V$ be the normalized volume of the simplex $L$. Then
$$\det M_L=\frac17(\det Q) V^2.$$
\end{lem}

\begin{proof}
We may assume that $m=0$ as translations do not change the inertia matrix or the normalized volume of $L$. Then $v_0=-(v_1+\ldots+v_6)$. Let $S$ be the $6\times 6$ matrix with columns $v_k$ for $k=1,\ldots,6$. Then for $k=1,\ldots, 6$, $v_k = S e_k$ where $e_k$ is the $k$th vector of the standard basis in $\R^6$ and $v_0 = -S\mathbf 1$ where $\mathbf 1 = (1,1,1,1,1,1)^t$. Combining this together, we get 
$$\sum_{k=0}^6 (v_k-m)(v_k-m)^t=\sum_{k=0}^6 v_kv_k^t= S\mathbf 1 \mathbf 1^t S^t+\sum_{k=1}^6 Se_ke_k^tS^t = S(I_6+\mathbf 1 \mathbf 1^t)S^t$$
where $I_6$ is the $6\times 6$ identity matrix.

Then 
\begin{multline*}
\det(M_L)=\det \left(A(\sum_{k=0}^6 (v_k-m)(v_k-m)^t)A^t \right)=\\=\det(A)(\det S)\det \left(I_6+ \mathbf 1 \mathbf 1^t \right)\det(S^t)\det(A^t).
\end{multline*}
Since $A^tA=Q$, $\det(I_6+ \mathbf 1 \mathbf 1^t)=7$, and $|\det S|=\frac 17V$, the lemma follows.
\end{proof}

As a corollary, $M_L$ is a positive definite matrix as it is positive semidefinite and has a positive determinant. Our next goal is to bound trace of the inertia matrix for simplices of the Delaunay triangulation $\mathcal D(Q)$.

\begin{lem}
\label{lem:trace-inertia}
Let $L$ be a simplex with circumradius at most $1$ with respect to the positive definite quadratic form $Q$. Then trace of the matrix $M_L$ is at most $7$.
\end{lem}
\begin{proof}
We use that $\tr(XY)=\tr(YX)$ for any appropriate matrices $X$ and $Y$. Then 
\begin{multline*}\tr(M_L)=\sum_{k=0}^6\tr\left(A(v_k-m)(v_k-m)^t A^t \right)=\sum_{k=0}^6\tr\left((v_k-m)^t A^t A(v_k-m)\right)=\\=\sum_{k=0}^6(v_k-m)^t Q(v_k-m).
\end{multline*}

Let $c$ be the circumcenter of $L$. Writing $v_k-c = (v_k-m)+(m-c)$ and applying the quadratic form $Q$ to both sides we obtain

\begin{multline*}
\sum_{k=0}^6(v_k-c)^t Q(v_k-c)=\\=\sum_{k=0}^6(v_k-m)^t Q(v_k-m)+7(m-c)^tQ(m-c)+2\left(\sum_{k=0}^6(v_k-m) \right)^tQ(m-c)=\\=\sum_{k=0}^6(v_k-m)^t Q(v_k-m)+7(m-c)^tQ(m-c).\end{multline*}

Therefore,
$$\tr(M_L)=\sum_{k=0}^6(v_k-m)^t Q(v_k-m)\leq \sum_{k=0}^6(v_k-c)^t Q(v_k-c) \leq 7$$ since for each $k$, $(v_k-c)^t Q(v_k-c)$ is the squared circumradius of $L$ and is at most $1$.
\end{proof}

We apply these properties of the trace and determinant to obtain lower bounds for the smallest eigenvalue of $M_L$ for simplices of the Delaunay triangulation of $Q$ provided $\det Q>4.395397$ and all simplices in $\mathcal D(Q)$ have circumradii at most $1$.

\begin{lem}
\label{lem:lambda}
Let $L$ be a simplex of the Delaunay triangulation of $Q$ with $\det Q>4.395397$. Suppose that the circumradius of $L$ is at most $1$. Let $\lambda$ be the smallest eigenvalue of $M_L$.
\begin{enumerate}
\item If normalized volume of $L$ is $2$, then $\lambda > 1.073933$.
\item If normalized volume of $L$ is $1$, then $\lambda > 0.128040$.
\end{enumerate}
\end{lem}

\begin{proof}
Since $\tr M_L\leq 7$, $\lambda \leq \frac 76$ and the sum of the remaining five eigenvalues of $M_L$ is at most $7-\lambda$. From the AM-GM inequality, 
$$ \det M_L \leq \lambda \left(\frac{7-\lambda}{5}\right)^5.$$
We note that the function $f(x)=x\left(\frac{7-x}{5}\right)^5$ is increasing on $\left[0,\frac76\right]$.

Now we turn our attention to our specific normalized volumes. Let $V$ be the normalized volume of $L$. If $V=2$, then from Lemma~\ref{lem:det-inertia}
$$\lambda \left(\frac{7-\lambda}{5}\right)^5 \geq \det M_L=\frac47\det Q > 2.511655.$$
From monotonicity of $f(x)$, we get $\lambda > 1.073933$.

If $V=1$, then from Lemma~\ref{lem:det-inertia}
$$\lambda \left(\frac{7-\lambda}{5}\right)^5 \geq \det M_L=\frac17\det Q > 0.627913.$$
From monotonicity of $f(x)$, we get $\lambda > 0.128040$.
\end{proof}

The next lemma connects length of a vector with respect to $Q$ with its representation as combination of the vertices of $L$ with sum zero coefficients. 

\begin{lem}
\label{lem:int-norm}
Let $L=\conv\{v_0,\ldots,v_6\}$ and let $\lambda$ be the smallest eigenvalue of $M_L$. Let $x=\sum_{k=0}^6\alpha_k v_k$ with $\sum_{k=0}^6\alpha_k=0$. Then $x^tQx\geq \lambda \sum_{k=0}^6 \alpha_k^2$.
\end{lem}

\begin{proof}
Let $\alpha = (\alpha_0,\ldots,\alpha_6)^t$. We first write $x$ as a combination of the vectors from the centroid $m=\frac17(v_0+\ldots+v_6)$ of $L$ to the vertices of $L$.
$$x=\sum_{k=0}^6\alpha_k v_k=\sum_{k=0}^6\alpha_k (v_k-m).$$

Then $x^tQx=\alpha^t G \alpha$ where $G$ is the $7\times 7$ matrix with entries $(v_k-m)^tQ(v_l-m)$ (or the Gram matrix of the vectors $v_k-m$ for $k=0,\ldots,6$ with respect to $Q$).

Let $A$ be the $6\times 6$ matrix such that $Q=A^tA$ and that was used to construct the inertia matrix $M_L$ in Definition~\ref{def:inertia}. Let $W$ be the $6\times 7$ matrix with columns $A(v_k-m)$ for $k=0,\ldots, 6$. Then the entries of $W^tW$ are $(A(v_k-m))^tA(v_l-m)=(v_k-m)^tA^tA(v_l-m)=(v_k-m)^tQ(v_l-m)$ so $W^tW=G$.  On the other hand, $WW^t=M_L$. 

This means that non-zero eigenvalues of $G$ and $M_L$ coincide and $\lambda$ is the smallest non-zero eigenvalue of $G$. The matrix $G$ has eigenvector $(1,1,1,1,1,1,1)^t$ with eigenvalue $0$; this eigenvector is orthogonal to $\alpha$ in the standard inner product. Decomposing $\alpha$ in a linear combination of an orthonormal basis of eigenvectors of $G$, we obtain $$x^tQx=\alpha^t G \alpha \geq \lambda \sum_{k=0}^6(\alpha_k^2).$$
\end{proof}

\subsection{Convex optimization}

Our last tool will be convex optimization for determinant of a positive definite quadratic form provided each simplex from a given set has circumradius at most $1$ with respect to that form. This approach was described by Sch\"urmann and Vallentin \cite{SV06} and here we mostly reproduce the details and adapt them to our computations in the package \cite{package}.

The following lemma is crucial for our computations. See \cite[Prop. 7.1]{SV06} for a proof.

\begin{lem}[Delaunay, Dolbilin, Ryshkov, Stogrin, \cite{DDRS70}]
\label{lem:B_L}
Let $L=\conv\{\mathbf 0,v_1,\ldots,v_6\}$ be a simplex in $\R^6$. Then the circumradius of $L$ with respect to the positive definite quadratic form $Q$ is at most $1$ if and only if the matrix
$$
B_L(Q):=\left(
\begin{array}{ccccc}
4 & v_1^tQv_1& v_2^tQv_2& \ldots & v_6^tQv_6\\
v_1^tQv_1& v_1^tQv_1& v_1^tQv_2& \ldots & v_1^tQv_6\\
v_2^tQv_2& v_1^tQv_2& v_2^tQv_2& \ldots & v_2^tQv_6\\
\vdots & \vdots & \vdots & \ddots & \vdots \\
v_6^tQv_6& v_1^tQv_6& v_2^tQv_6& \ldots & v_6^tQv_6
\end{array}
\right)
$$
is positive semidefinite.
\end{lem}

For a fixed simplex $L$, the matrix $B_L(Q)$ is a linear function in the entries of $Q$. That is, there are matrices $G_0$ and $G_{ij}$ for $1\leq i\leq j \leq 6$ such that 
$$B_L(Q)=G_0+\sum_{1\leq i\leq j \leq 6}q_{ij}G_{ij}$$
where $q_{ij}$ are the entries of $Q$. Therefore, for a fixed simplex $L$, the collection $\mathcal S_L$ of all positive definite quadratic forms $Q$ such that $B_L(Q)$ is positive semidefinite is a convex set in $\mathcal S_6^{>0}$. 

If $\mathcal L$ is a finite collection of simplices, then the family of all positive definite quadratic forms $Q$ such that all simplices in $\mathcal L$ have circumradii at most $1$ is convex as well as an intersection of convex sets. This allows us to maximize the concave function $\log \det Q$  (which is concave on $\mathcal S_6^{>0}$) over the convex set $\bigcap_{L\in \mathcal L}\mathcal S_L$. We also note that if $\mathcal L$ is the set of all simplices from a Delaunay triangulation of $\Z^6$, then this approach can be used to optimize the lattice sphere covering over the closure of the corresponding secondary cone by maximizing $\log \det Q$ over the convex set defined by the closure of the secondary cone and the intersection of the sets $\mathcal S_L$.

The convex optimization in \cite{package} is done using {\tt python}'s {\tt cvxpy} package. This optimization gives an approximate {\tt float} bound and we can use the output of the {\tt cvxpy} solver to search for a certified rational bound similarly to the proof of Theorem~\ref{thm:SV} in \cite{SV06}. Theoretical details about the certification process can be found in \cite[Sec. 6]{SV06} and here we describe an approach to give an upper bound for the determinant of $Q$; we illustrate how the process works for two simplices in rational arithmetic in Example~\ref{exam:two} below.

Let 
$$N=\left(
\begin{array}{cccc}
n_{00}& n_{01}& \ldots &n_{06} \\
n_{01}& n_{11}& \ldots &n_{16} \\
\vdots& \vdots& \ddots &\vdots \\
n_{06}& n_{16}& \ldots &n_{66} 
\end{array}    
\right)$$
be any positive $7\times 7$ semidefinite matrix. Let $D=\diag(n_{01},\ldots,n_{06})$ be the diagonal matrix with entries $n_{01},\ldots,n_{06}$ on the diagonal and $N'=(n_{ij})_{i,j=1}^6$ be the $6\times 6$ matrix obtained from $N$ by removing the first row and the first column. 

As in Lemma~\ref{lem:B_L}, let $L=\conv\{\mathbf 0,v_1,\ldots,v_6\}$. If the matrix $B_L(Q)$ is positive semidefinite and the circumradius of $L$ is at most $1$, then we may use the matrix $N$ to obtain an upper bound on $\det Q$ as follows. Let $V$ be the $6\times 6$ matrix with columns $v_1,\ldots,v_6$.

Since both $N$ and $B_L(Q)$ are positive semidefinite, 
$$0\leq \trace (NB_L(Q)) = 4n_{00} + 2\sum_{i=1}^6n_{0i}v_i^tQv_i+\trace(N'V^tQV).$$

Rewriting $\sum_{i=1}^6n_{0i}v_i^tQv_i=\trace(DV^tQV)=\trace(VDV^tQ)$ and $\trace(N'V^tQV)= \trace(VN'V^tQ)$, we obtain
\begin{equation}
\label{eq:trace}
-\trace(V(N'+2D)V^tQ)\leq 4n_{00}.
\end{equation}

Let $P_{L,N}=-V(N'+2D)V^t$. Note that this is a $6 \times 6$ symmetric matrix that depends on the simplex $L$ (as the matrix $V$ depends on $L$) and the matrix $N$.

Now suppose that we want to get an upper bound for $\det Q$ such that all matrices $B_L(Q)$ for $L\in \mathcal L$ are positive semidefinite where $\mathcal L$ is a finite collection of simplices\footnote{We do not restrict $\mathcal L$ to subsets of simplices from one Delaunay triangulation. They can even overlap as long as each simplex has $(0,0,0,0,0,0)^t$ as a vertex.}. Let $\{N_L:L\in \mathcal L\}$ be a collection of positive semidefinite $7\times 7$ matrices. Let ${\displaystyle C_{\mathcal L}=4\sum_{N_L}n_{00}}$ be the sum of all right-hand sides from \eqref{eq:trace} over all simplices $L\in\mathcal L$ and ${\displaystyle P_{\mathcal L}=\sum_{L\in \mathcal L}P_{L,N_L}}$. Then the sum of inequalities \eqref{eq:trace} of all $L\in \mathcal L$ gives us the inequality
$$\trace(P_{\mathcal L}Q)\leq C_{\mathcal L}.$$

If $C_{\mathcal L}> 0 $ and if $P_{\mathcal L}$ is positive definite, then we can get an upper bound for $\det Q$. Writing $A=Q^{\frac12}$, we get $$\trace(AP_{\mathcal L}A)=\trace(P_{\mathcal L}Q)\leq C_{\mathcal L}.$$
The matrix $AP_{\mathcal L}A$ is positive definite, so its eigenvalues are positive real numbers. From the AM-GM inequality,
$$\det Q \cdot \det P_{\mathcal L} = \det(AP_{\mathcal L}A)\leq \left(\frac{C_{\mathcal L}}{6}\right)^6 \qquad \text{or} \qquad \det Q \leq \frac{C_{\mathcal L}^6}{6^6 \det P_{\mathcal L}}.$$
So for a given collection of simplices $\mathcal L$, in order to get a good upper bound for $\det Q$, we need to choose an appropriate family of matrices $N_L,L\in \mathcal L$. We can use approximate solution for the dual convex optimization problem and this dual solution is also obtained through the {\tt cvxpy} package.

We illustrate how our computational approach works for a pair of two unimodular simplices and show that these two simplices force a small determinant for any relevant positive definite quadratic form.

\begin{exam}
\label{exam:two}
Let $\sigma$ be the simplex with vertices $e_1,\ldots,e_6$ and $x=(2,2,1,-1,-1,-1)$. We will show that if Delaunay triangulation of $Q$ contains both the standard simplex $\sigma_0$ with vertices at the origin and at all six basis vectors and $\sigma$, then $\det Q< 4.395397$ and therefore $Q$ gives worse lattice sphere covering density than $Q_6^{c}$.

To do that, we solve the convex optimization problem of maximizing $\log \det Q$ under constraints that both $B_{\sigma_0}$ and $B_\sigma$ are positive semidefinite as both simplices must have circumradii at most $1$; we note that we need to translate one of the vertices of $\sigma$ to the origin to use Lemma~\ref{lem:B_L}. This convex optimization problem is solved using {\tt python}'s {\tt cvxpy} package in the {\tt maxdet} function in the {\tt common/relax.py} file in package \cite{package}. The function's output contains approximate maximum $4.065983080565106$ of the determinant\footnote{This maximum and other values obtained through computations in {\tt float} arithmetic may vary slightly depending on the order of vertices we use for each simplex and order of simplices in the input of the {\tt maxdet} function.}, two approximate matrices 
{
\tiny 
\begin{equation*} N_{\sigma_0}= \left(
\begin{array}{rrrrrrr}
 0.74996228& -0.15318192& -0.15317792& -0.21112281& -0.21115639& -0.21115643& -0.21115359\\
-0.15318192&  0.03128785&  0.03128703&  0.04312243&  0.04312929& 0.0431293 &  0.04312872\\
-0.15317792&  0.03128703&  0.03128621&  0.04312131&  0.04312816& 0.04312817&  0.04312759\\
-0.21112281&  0.04312243&  0.04312131&  0.05943344&  0.0594429 & 0.05944291&  0.05944211\\
-0.21115639&  0.04312929&  0.04312816&  0.0594429 &  0.05945235& 0.05945236&  0.05945156\\
-0.21115643&  0.0431293 &  0.04312817&  0.05944291&  0.05945236& 0.05945238&  0.05945158\\
-0.21115359&  0.04312872&  0.04312759&  0.05944211&  0.05945156& 0.05945158&  0.05945078\\
\end{array}
\right)
\end{equation*}}
and 
{\tiny 
\begin{equation*}
N_\sigma=\left(
\begin{array}{rrrrrrr}
 0.75003782& -0.15323402& -0.15323811& -0.21117639& -0.21114801& -0.2111486& -0.21114181\\
-0.15323402&  0.03130598&  0.03130682&  0.0431437 &  0.0431379 & 0.04313802&  0.04313664\\
-0.15323811&  0.03130682&  0.03130765&  0.04314485&  0.04313905& 0.04313917&  0.04313778\\
-0.21117639&  0.0431437 &  0.04314485&  0.05945763&  0.05944963& 0.0594498 &  0.05944789\\
-0.21114801&  0.0431379 &  0.04313905&  0.05944963&  0.05944165& 0.05944181&  0.0594399 \\
-0.2111486 &  0.04313802&  0.04313917&  0.0594498 &  0.05944181& 0.05944198&  0.05944007\\
-0.21114181&  0.04313664&  0.04313778&  0.05944789&  0.0594399 & 0.05944007&  0.05943816\\
\end{array}    
\right)
\end{equation*}}
of $7\times 7$ Lagrange multipliers (optimizers for the dual convex programming problem) corresponding to $\sigma_0$ and $\sigma$, respectively, and an approximate matrix $Q$ that maximizes the determinant with $\det Q \approx 4.065983080565106$. Note that the two matrices $N_{\sigma_0}$ and $N_\sigma$ are almost equal; for the exact solution with these two simplices, they will be equal due to a symmetry of the configuration of two simplices.

Since the solution is approximate, we separately certify that the maximum of the determinant is less than $4.395397$ using rational arithmetic in the function {\tt certificate_bound} in {\tt common/certificate.py} in \cite{package} as follows. As output of {\tt maxdet} function, each matrix $N$ has {\tt float} entries. We transform each matrix $N$ to a matrix with rational entries that is guaranteed to be positive semidefinite in the function {\tt rational_psd} in {\tt common/relax.py} in \cite{package} as follows.
\begin{itemize}
    \item We replace $N$ with $N^o=0.5(N+N^t)$ which is a symmetric matrix. While $N$ is expected to be symmetric, the approximate values may make it non-symmetric due to approximation error.
    \item Find eigenvalues $\lambda_1,\ldots,\lambda_7$ of $N^o$ and (theoretically) orthogonal basis of eigenvectors $v_1,\ldots,v_7$ of $N^o$ numerically. Note that some of the eigenvalues may be negative due to an approximation error.
    \item We replace each negative $\lambda_i$ with $\lambda_i'=0$ and approximate the remaining eigenvalues $\lambda_j$ with rational numbers $\lambda_j'$ with a common denominator $10^2$ (default denominator for our computations is $10^6$ and we have chosen a smaller denominator for this example to get matrices that are easier to read). We replace each eigenvector $v_i$ with $v_i'$ for $i=1,\ldots, 7$ by approximating each coordinate with a rational number with the same denominator $10^2$.
    \item The matrix $N^o$ is then replaced with $$N^{\Q}=\sum_{i=1}^7\lambda_i' v_i' (v_i')^t$$
    which is a positive semidefinite matrix with rational entries approximating both $N$ and $N^o$.
\end{itemize}
After performing these steps we get two positive semidefinite matrices $N_{\sigma_0}^\Q$ and $N_\sigma^\Q$ with rational entries. While these two matrices may not be the precise output of {\tt maxdet}, we may use them for the certification process as described before this example. For the two matrices $N_{\sigma_0}$ and $N_\sigma$ above, the resulting rational matrices are

\begin{equation*} 
\renewcommand{\arraystretch}{1.5} 
N^{\mathbb Q}_{\sigma_0}= \left(
\begin{array}{rrrrrrr}
\frac{6069}{8000} & \frac{-6069}{40000} & \frac{-6069}{40000} & \frac{-1071}{5000} & \frac{-1071}{5000} & \frac{-1071}{5000} & \frac{-1071}{5000} \\ 
\frac{-6069}{40000} & \frac{6069}{200000} & \frac{6069}{200000} & \frac{1071}{25000} & \frac{1071}{25000} & \frac{1071}{25000} & \frac{1071}{25000} \\ 
\frac{-6069}{40000} & \frac{6069}{200000} & \frac{6069}{200000} & \frac{1071}{25000} & \frac{1071}{25000} & \frac{1071}{25000} & \frac{1071}{25000} \\ 
\frac{-1071}{5000} & \frac{1071}{25000} & \frac{1071}{25000} & \frac{189}{3125} & \frac{189}{3125} & \frac{189}{3125} & \frac{189}{3125} \\ 
\frac{-1071}{5000} & \frac{1071}{25000} & \frac{1071}{25000} & \frac{189}{3125} & \frac{189}{3125} & \frac{189}{3125} & \frac{189}{3125} \\ 
\frac{-1071}{5000} & \frac{1071}{25000} & \frac{1071}{25000} & \frac{189}{3125} & \frac{189}{3125} & \frac{189}{3125} & \frac{189}{3125} \\ 
\frac{-1071}{5000} & \frac{1071}{25000} & \frac{1071}{25000} & \frac{189}{3125} & \frac{189}{3125} & \frac{189}{3125} & \frac{189}{3125}
\end{array}
\right)
\end{equation*}
and
\begin{equation*} 
\renewcommand{\arraystretch}{1.5} 
N^{\mathbb Q}_{\sigma}= \left(
\begin{array}{rrrrrrr}
\frac{6069}{8000} & \frac{-6069}{40000} & \frac{-6069}{40000} & \frac{-1071}{5000} & \frac{-1071}{5000} & \frac{-1071}{5000} & \frac{-1071}{5000} \\ 
\frac{-6069}{40000} & \frac{6069}{200000} & \frac{6069}{200000} & \frac{1071}{25000} & \frac{1071}{25000} & \frac{1071}{25000} & \frac{1071}{25000} \\ 
\frac{-6069}{40000} & \frac{6069}{200000} & \frac{6069}{200000} & \frac{1071}{25000} & \frac{1071}{25000} & \frac{1071}{25000} & \frac{1071}{25000} \\ 
\frac{-1071}{5000} & \frac{1071}{25000} & \frac{1071}{25000} & \frac{189}{3125} & \frac{189}{3125} & \frac{189}{3125} & \frac{189}{3125} \\ 
\frac{-1071}{5000} & \frac{1071}{25000} & \frac{1071}{25000} & \frac{189}{3125} & \frac{189}{3125} & \frac{189}{3125} & \frac{189}{3125} \\ 
\frac{-1071}{5000} & \frac{1071}{25000} & \frac{1071}{25000} & \frac{189}{3125} & \frac{189}{3125} & \frac{189}{3125} & \frac{189}{3125} \\ 
\frac{-1071}{5000} & \frac{1071}{25000} & \frac{1071}{25000} & \frac{189}{3125} & \frac{189}{3125} & \frac{189}{3125} & \frac{189}{3125}
\end{array}
\right).
\end{equation*}

Using the two matrices $N_{\sigma_0}^{\mathbb Q}$ and $N_{\sigma}^{\mathbb Q}$ in the process described before this example, we obtain the following rational upper bound for $\det Q$.
$$\det Q \leq \frac{6411541765625}{1576312897536} \approx 4.067429617334951.$$
While the bound is slightly greater than the approximate solution of the convex optimization problem, it allows us to exclude facets of the Type $(2,2,1,-1,-1,-1)$ from potential Delaunay triangulations of quadratic forms that have covering density below $2.4649$.
\end{exam}

We note that we accept the rational bound for the determinant only if it is below $4.395397$. If it is above the threshold, the constant $C_{\mathcal L}$ is not positive, or the matrix $P_{\mathcal L}$ is not positive definite, then we do not get a satisfactory upper bound and treat the collection of simplices $\mathcal L$ as theoretically possible in a Delaunay triangulation of $Q$ with covering density below $2.4649$.

In further sections we use these tools to provide details of the specific steps from the proof of Theorem~\ref{thm:main}.

\section{Volume bound}
\label{sec:volume}

In this section we prove that if $Q$ is a generic six-dimensional positive definite quadratic form such that $\Theta(Q)<2.4649$, then the Delaunay triangulation of $Q$ is unimodular. That is, all simplices of $\mathcal D(Q)$ have normalized volume $1$. Recall that we assume that the covering radius of $Q$ is $1$ and this implies that circumradii of all simplices in $\mathcal D(Q)$ are at most $1$ with respect to $Q$.

We start from a preliminary bound.

\begin{lem}
\label{lem:vol2}
Normalized volumes of all simplices in the Delaunay triangulation of $Q$ do not exceed $2$.
\end{lem}

\begin{proof}
Let $L$ be a simplex of normalized volume $V$ in the Delaunay triangulation of $Q$; $V$ is a positive integer. Let $R$ be the circumradius of $L$ with respect to $Q$. We map the ellipsoid circumscribed around $L$ into a Euclidean ball of radius $R$ by an affine transformation. This turns $L$ into a simplex of Euclidean volume  $\frac{V\cdot \sqrt{\det Q}}{6!}$ inscribed in a sphere of radius $R$. This volume does not exceed the volume of the regular simplex with circumradius $1$ and therefore
$$\frac{V\cdot \sqrt{\det Q}}{6!}\leq \frac{\sqrt 7}{6!}\cdot \left(\frac{7}{6}\right)^3.$$

Thus, $V\leq \frac{343\sqrt 7}{216\sqrt{\det Q}}$ and since the covering radius $R(Q)= 1$, we get $\frac{\kappa_6}{\sqrt{\det Q}}= \Theta(Q)<2.5$. Combining the inequalities and $\kappa_6=\frac{\pi^3}{6}$, we obtain
$$V \leq \frac{1715\sqrt{7}}{72\pi^3}\approx 2.0325 < 3.$$
\end{proof}

\begin{lem}
\label{lem:vol1}
The Delaunay triangulation of $Q$ is unimodular.
\end{lem}

\begin{proof}
From Lemma~\ref{lem:vol2}, it is enough to show that no simplex of the Delaunay triangulation $\mathcal D(Q)$ has normalized volume $2$. Assume there is such a simplex $L$ in $\mathcal D(Q)$. 

We assume that one of the vertices of $L$ is the origin and that the remaining $6$ vertices of $L$ form a $6\times 6$ matrix in Hermite Normal Form . In other words, we may assume that six of the vertices of $L$ are 
\begin{center}
\begin{tabular}{ll}
$O=(0,0,0,0,0,0)^t$, & $A=(1,0,0,0,0,0)^t$, \\     
$B=(0,1,0,0,0,0)^t$, & $C=(0,0,1,0,0,0)^t$, \\     
$D=(0,0,0,1,0,0)^t$, & $E=(0,0,0,0,1,0)^t$, \\     
\end{tabular}
\end{center}
and the final vertex is $F=(*,*,*,*,*,2)^t$ where each $*$ can be either $0$ or $1$ independently from the others.

We consider six cases based on the number of zeros among the coordinates of $F$. It is enough to consider only one representative for each case as all others can be obtained by permutations of coordinates.

{\bf Case 1.} If $F=(0,0,0,0,0,2)^t$, then the midpoint of $OF$ is a lattice point, so $OF$ cannot be an edge in any Delaunay triangulation of $\Z^6$.

{\bf Case 2.} If $F=(1,0,0,0,0,2)^t$, then the midpoint of $AF$ is a lattice point, so $AF$ cannot be an edge in any Delaunay triangulation of $\Z^6$.

{\bf Case 3.} If $F=(1,1,0,0,0,2)^t$, then the parity classes of edges $OF$ and $AB$ coincide but the edges are not parallel. This is not possible in Delaunay triangulation of $\Z^6$ due to Lemma~\ref{lem:parity}.

{\bf Case 4.} If $F=(1,1,1,0,0,2)^t$, then the parity classes of edges $AF$ and $BC$ coincide but the edges are not parallel. This is not possible in Delaunay triangulation of $\Z^6$ due to Lemma~\ref{lem:parity}.

{\bf Case 5.} If $F=(1,1,1,1,0,2)^t$, then we note that the affine transformation
$$
\left(
\begin{array}{r}x_1\\x_2\\x_3\\x_4\\x_5\\x_6    
\end{array}
\right)\mapsto
\left(
\begin{array}{rrrrrr}
1 & 0& 0& 0& 0& 0\\
0 & 1& 0& 0& 0& 0\\
0 & 0& 1& 0& 0& 0\\
0 & 0& 0& 1& 0& 0\\
-1&-1&-1&-1&-1& 2\\
0 & 0& 0& 0& 0& 1
\end{array}
\right)
\left(
\begin{array}{r}x_1\\x_2\\x_3\\x_4\\x_5\\x_6    
\end{array}
\right)+\left(
\begin{array}{r}0\\0\\0\\0\\1\\0
\end{array}
\right)
$$
(with determinant $-1$) transforms this simplex into the one from {\bf Case 6} below while giving bijection of the lattice $\Z^6$ with itself. So we just consider the remaining case as the formulas are slightly more symmetric.

{\bf Case 6.} If $F=(1,1,1,1,1,2)^t$, let $P\neq O$ be a point in $\Z^6$ such that the simplex $L'=PABCDEF$ is also in the Delaunay triangulation of $Q$. Let the barycentric coordinates of $P$ with respect to the simplex $L$ be $(m_0,\ldots,m_6)$.

Since the circumradius of $L'$ is at most $1$, the lengths of vectors $PA,PB,PC,PD,PE,PF$ with respect to the quadratic form $Q$ are at most $2$. Let $v=PA$, then from the barycentric coordinates of $P$ we get
$$v=m_0O+(m_1-1)A+m_2B+m_3C+m_4D+m_5E+m_6F$$
and since the sum of coefficients is 0, we use Lemmas~\ref{lem:int-norm} and~\ref{lem:lambda} to get
$$m_0^2+(m_1-1)^2+m_2^2+m_3^2+m_4^2+m_5^2+m_6^2\leq \frac{v^tQv}{\lambda}<\frac{4}{1.073933}<\frac{15}{4}$$
where $\lambda$ is the smallest eigenvalue of the inertia matrix $M_L$. This can be rewritten as 
$$m_0^2+m_1^2+m_2^2+m_3^2+m_4^2+m_5^2+m_6^2<\frac{11}{4}+2m_1.$$

Applying the same arguments to vectors $PB,PC,PD,PE$, and $PF$ we obtain that for $k=1,\ldots,6$,

$$m_0^2+m_1^2+m_2^2+m_3^2+m_4^2+m_5^2+m_6^2<\frac{11}{4}+2m_k.$$

Now we write the actual coordinates of $P$ from its barycentric coordinates
$$P=(m_1+m_6,m_2+m_6,m_3+m_6,m_4+m_6,m_5+m_6,2m_6)^t \quad \text{and} \quad \sum_{k=0}^6m_k=1.$$
Note that all $m_k$ are integers or half integers because the normalized volume of $L$ is $2$ and $m_0<0$ because $P$ and $O$ have to be on different sides of the hyperplane $ABCDEF$. All $m_k$ for $k=1,\ldots,6$ have to be in $\Z$ or all in $\Z+\frac12$ simultaneously to make coordinates of $P$ into integers. This means $m_0$ is an integer and $m_0\leq -1$. 

We consider two cases. If all $m_k$ for $k=1,\ldots,6$ are integers, then $m_0^2+m_1^2\geq 2m_1$ because $m_0^2\geq 1$. So $m_2^2+\ldots+m_6^2<3$ and at least one of these $m_k$ is $0$. This implies that 
$m_0^2+\ldots+m_6^2<3$ but $m_0^2\geq 1$ and sum of the other $m_k$'s is at least $2$ so sum of their squares is at least $2$ because all are integers. This is a contradiction.

Now suppose all $m_k$ for $k=1,\ldots,6$ are half-integers. If all these $m_k$ are positive, then $m_0= 1 - (m_1+\ldots+m_6) \leq -2$ and 
$$m_0^2+(m_1-1)^2+m_2^2+m_3^2+m_4^2+m_5^2+m_6^2 \geq 4 +6\cdot \frac14>\frac{15}{4}$$
which is a contradiction. If at least one of $m_k$ for $k=1,\ldots,6$ is negative, then for this $k$, $2m_k\leq -1$ and
$$m_0^2+m_1^2+m_2^2+m_3^2+m_4^2+m_5^2+m_6^2\geq 1+6\cdot \frac 14 = \frac{10}{4} > \frac{11}{4} + 2m_k$$ which is again a contradiction.
\end{proof}

So from now on we may assume that the Delaunay triangulation of $Q$ is unimodular if  $\Theta(Q)<2.4649$.

\section{Possible types of neighbors}
\label{sec:neighbors}

In this section we prove that there are three possible types of facets in the Delaunay triangulation for the form $Q$. These are Type A, Type B, and Type C from Section~\ref{sec:lattice}. We note that we provide more computational details and outputs than in the further sections to make sure that the reader is familiar with the process of our computations. Further sections concentrate on the logic of the computations while details of the outputs are skipped in the paper but supplied in the package \cite{package}.

The process is organized as follows. Recall that we may assume that the Delaunay triangulation of $Q$ is unimodular. We also may assume that one of the simplices of the triangulation is the standard simplex $$\sigma_0 = \conv \left(\begin{array}{l}
(0,0,0,0,0,0)^t,\\
(1,0,0,0,0,0)^t,\\
(0,1,0,0,0,0)^t,\\
(0,0,1,0,0,0)^t,\\
(0,0,0,1,0,0)^t,\\
(0,0,0,0,1,0)^t,\\
(0,0,0,0,0,1)^t \\
\end{array}\right)$$
with vertices given by the origin $(0,0,0,0,0,0)^t$ and the standard basis vectors $e_1,\ldots,e_6$. Using Lemma~\ref{lem:int-norm}, we may limit potential edges in the Delaunay triangulation of $Q$ and then limit potential neighbors of $\sigma_0$. Finally, for the potential neighbors of $\sigma_0$, we use convex optimization approach described in Section~\ref{sec:lemmas} to discard local configurations that guarantee small determinant of $Q$.

\begin{lem}
\label{lem:edges}
Let $e$ be an edge of the Delaunay triangulation of $Q$. Then $e$, up to change of sign, belongs to a list of $17822$ edges.
\end{lem}
\begin{proof}
Let $x=(x_1,\ldots,x_6)^t$ be the vector representing $e$. We apply Lemma~\ref{lem:int-norm} to $x$ and simplex $\sigma_0$. We can write 
$$x=-(x_1+\ldots+x_6)\mathbf 0 + x_1e_1+\ldots+x_6e_6$$
where $\mathbf 0$ is the zero vector. This is exactly the representation needed for Lemma~\ref{lem:int-norm}. Since $x$ represents an edge in the Delaunay triangulation of $Q$ with covering radius $1$, $x^tQx\leq 4$ and Lemma~\ref{lem:int-norm} gives us the bound 
$$(x_1+\ldots+x_6)^2+x_1^2+\ldots+x_6^2 \leq \frac{4}{\lambda}$$
where $\lambda$ is the smallest eigenvalue of the inertia matrix $M_{\sigma_0}$.

From Lemma~\ref{lem:lambda} we have $\lambda > 0.128040$ and so
$$(x_1+\ldots+x_6)^2+x_1^2+\ldots+x_6^2 \leq \frac{4}{\lambda}<32.$$
In particular, for each $k=1,\ldots,6$, $|x_k|\leq 5$.

This gives finitely many options to go through. Once we remove all candidates where all $x_k, k=1,\ldots,6$ are even and the candidates from the parity classes of the edges of $\sigma_0$ except the edges of $\sigma_0$ themselves, we get the list of $17822$ edges.
This is done in the accompanying package in the script {\tt common/edges.py} \cite{package}.
\end{proof}

\begin{lem}
\label{lem:23apexes}
Let $\sigma$ be the neighbor of $\sigma_0$ in the Delaunay triangulation of $Q$ sharing the facet with vertices $e_1,\ldots,e_6$. Let $x$ be the last vertex of $\sigma$. Then $x$ belongs to the following list of $23$ points
{\small
$$\begin{array}{llll}
(1, 1, 0, 0, 0, 0)^t,&
(1, 1, 1, 0, 0, -1)^t,&
(1, 1, 1, 1, -1, -1)^t,&
(1, 1, 1, 1, 0, -2)^t,\\
(1, 1, 1, 1, 1, -3)^t,&
(2, 1, 1, 0, -1, -1)^t,&
(2, 1, 1, 1, -1, -2)^t,&
(2, 1, 1, 1, 0, -3)^t,\\
(2, 2, 1, -1, -1, -1)^t,&
(2, 2, 1, 1, -1, -3)^t,&
(3, 1, 0, 0, -1, -1)^t,&
(3, 1, 1, -1, -1, -1)^t,\\
(3, 1, 1, 0, -1, -2)^t,&
(3, 1, 1, 0, 0, -3)^t,&
(3, 1, 1, 1, -2, -2)^t,&
(3, 1, 1, 1, -1, -3)^t,\\
(3, 2, 0, -1, -1, -1)^t,&
(3, 2, 1, -1, -1, -2)^t,&
(3, 3, -1, -1, -1, -1)^t,&
(3, 3, 0, -1, -1, -2)^t,\\
(4, 1, 0, -1, -1, -1)^t,&
(4, 1, 1, -1, -1, -2)^t,&
(4, 2, -1, -1, -1, -1)^t&
\end{array}$$}
after possible permutation of coordinates.
\end{lem}

\begin{proof}
Since all simplices in the Delaunay triangulation of $Q$ are unimodular, the sum of all coordinates of $x$ is $2$. For $e_1$, we consider all $17822$ potential edges from Lemma~\ref{lem:edges} and their negatives starting from $e_1$ and then check whether the endpoint has sum of coordinates $2$ and whether the other five edges connecting the endpoint with $e_2,\ldots,e_6$ are in the list from Lemma~\ref{lem:edges}.

Endpoints satisfying these conditions are given in {\tt 02_step2_types/step2_apexes.json} which can be obtained running script {\tt 02_step2_types/01_apex.py} in  \cite{package}.
\end{proof}

This list restricts potential types of facets in the Delaunay triangulation $\mathcal D(Q)$ to just 23 types as we can apply an affine transformation to change any given unimodular simplex to $\sigma_0$ and its selected facet to the facet of $\sigma_0$ opposite to the origin while preserving $\Z^6$. Our next goal is to discard almost all potential types of facets except a list of three Types A, B, and C. We do that following the process of Example~\ref{exam:two} which helps eliminating all but five types first.

\begin{lem}
\label{lem:5types}
Under the conditions of Lemma~\ref{lem:23apexes}, if the common facet of $\sigma_0$ and $\sigma$ does not belong to the following five types, then $\det Q < 4.395397$.
\begin{center}
\begin{tabular}{ll}
Type A:& $(1,1,1,0,0,-1)$\\
Type B:& $(1,1,0,0,0,0)$\\
Type C:& $(1,1,1,1,-1,-1)$\\
Type D:& $(2,1,1,0,-1,-1)$\\
Type E:& $(1,1,1,1,0,-2)$
\end{tabular}
\end{center}
\end{lem}
\begin{proof}
The only tool we need to use here is the convex optimization bound for the pair of simplices $(\sigma_0,\sigma)$ from Lemma~\ref{lem:23apexes}; an example of how this bound is used is given in Example~\ref{exam:two}.

We run this approach over all 23 potential simplices $\sigma$ representing different facet types in the script {\tt 02_step2_types/02_types_ABCDE.py} in \cite{package}. The output of the script is the following.

{\tiny\tt\begin{center}
\begin{tabular}{l}
[1, 1, 0, 0, 0, 0] - Accepted. Solver status = optimal. Convex optimization bound= 13.013247861218389\\ \relax
[1, 1, 1, 0, 0, -1] - Accepted. Solver status = optimal. Convex optimization bound= 9.377136974266774\\ \relax
[1, 1, 1, 1, -1, -1] - Accepted. Solver status = optimal. Convex optimization bound= 7.307362067600618\\ \relax
[1, 1, 1, 1, 0, -2] - Accepted. Solver status = optimal. Convex optimization bound= 5.40426734454741\\ \relax
[1, 1, 1, 1, 1, -3] - Rejected. Solver status = optimal. Certified bound= 2.6117148674158597\\ \relax
[2, 1, 1, 0, -1, -1] - Accepted. Solver status = optimal. Convex optimization bound= 5.4042673720529555\\ \relax
[2, 1, 1, 1, -1, -2] - Rejected. Solver status = optimal. Certified bound= 4.065983120088399\\ \relax
[2, 1, 1, 1, 0, -3] - Rejected. Solver status = optimal. Certified bound= 2.6728812991114075\\ \relax
[2, 2, 1, -1, -1, -1] - Rejected. Solver status = optimal. Certified bound= 4.065983118032878\\  \relax
[2, 2, 1, 1, -1, -3] - Rejected. Solver status = optimal. Certified bound= 2.4611726791723925\\  \relax
[3, 1, 0, 0, -1, -1] - Rejected. Solver status = optimal. Certified bound= 2.1569433989759137\\ \relax
[3, 1, 1, -1, -1, -1] - Rejected. Solver status = optimal. Certified bound= 2.6117148688775145\\ \relax
[3, 1, 1, 0, -1, -2] - Rejected. Solver status = optimal. Certified bound= 2.6728813017611084\\ \relax
[3, 1, 1, 0, 0, -3] - Rejected. Solver status = optimal. Certified bound= 2.365358637892606\\ \relax
[3, 1, 1, 1, -2, -2] - Rejected. Solver status = optimal. Certified bound= 2.461172680931365\\ \relax
[3, 1, 1, 1, -1, -3] - Rejected. Solver status = optimal. Certified bound= 2.238683746278662\\ \relax
[3, 2, 0, -1, -1, -1] - Rejected. Solver status = optimal. Certified bound= 2.672881298564234\\ \relax
[3, 2, 1, -1, -1, -2] - Rejected. Solver status = optimal_inaccurate. Certified bound= 2.4611727087100443\\ \relax
[3, 3, -1, -1, -1, -1] - Rejected. Solver status = optimal_inaccurate. Certified bound= 2.238683760275946\\ \relax
[3, 3, 0, -1, -1, -2] - Rejected. Solver status = optimal. Certified bound= 2.1102911161224247\\ \relax
[4, 1, 0, -1, -1, -1] - Rejected. Solver status = optimal. Certified bound= 0.6415586132770048\\ \relax
[4, 1, 1, -1, -1, -2] - Rejected. Solver status = optimal. Certified bound= 1.300405271047154\\ \relax
[4, 2, -1, -1, -1, -1] - Rejected. Solver status = optimal. Certified bound= 1.3004052693305161\\ \relax
Total remaining 5 candidates  
\end{tabular}
\end{center}}

We note that even when the {\tt cvxpy} solver status was {\tt optimal_inaccurate}, we were able to find a pair of positive semidefinite matrices that certify the bound for the determinant. We also note that the upper bound we obtained for the type $(2,2,1,-1,-1,-1)$ is better than the bound from Example~\ref{exam:two}; this is because we use the default denominator $10^6$ to turn {\tt float} matrices into rational positive semidefinite matrices in the script versus the $10^2$ denominator in Example~\ref{exam:two}.
\end{proof}

This lemma also gives a list of potential neighbors for every simplex $S$ in the Delaunay triangulation of $Q$. We just need to map $S$ to the standard simplex by an affine transformation or, equivalently, consider barycentric coordinates of the additional vertex of the potential neighbor of $S$.

In order to exclude two additional facet Types, Type D and Type E, we need to use a deeper search among potential configurations of simplices. Let us start with a description of the configurations of simplices we use and what we expect from these configurations.

\begin{defin}\label{def:conf}
Let $S$ be a unimodular simplex with vertices in $\Z^6$. The {\it equivalence class} of $S$ (also denoted $S$) is the collection of all simplices obtained from $S$ by integer translations and central symmetries at half-integer points.

We call a finite collection $\{S_1,\ldots,S_n\}$ of equivalence classes of simplices a {\it configuration}, if it satisfies the following three properties:
\begin{enumerate}[label=(\alph*)]
\item Each non-trivial parity class of lattice vectors contains at most one representative among edges of $S_1,\ldots,S_n$; two opposite vectors are considered as the same representative. And no edge from any simplex in the configuration belongs to the trivial parity class.
\item Each facet $F$ of each simplex from $S_i$ belongs to at most two simplices from all classes in the configuration, and if it belongs to two simplices, then these simplices are on opposite sides of $F$.
\item For every $i$, if we map a simplex from $S_i$ to the standard simplex $\sigma_0$ by an affine transformation, then the edges of all other simplices in the configuration belong to the list of 17822 edges from Lemma~\ref{lem:edges}. This condition will be called the {\it frame condition}.
\end{enumerate}
\end{defin}
    
We note that if $Q$ is our generic quadratic form with $\det Q> 4.395397$ and covering radius at most $1$, then these three conditions are necessary for the Delaunay triangulation of $Q$. Also, the configuration given by the Delaunay triangulation must contain exactly $360$ different classes of simplices. Indeed, there should be $6!=720$ different translational classes simplices of $\mathcal D(Q)$ from the volume argument, and half of that once we incorporate central symmetries.

We do not require any stronger properties for non-intersection of simplices even locally as it is more complicated to check computationally. For example, we do not track whether two simplices in one configuration intersect by an internal point of both unless they share a facet; in that case, the two simplices must be on different sides of the facet.

\begin{defin}
Let $\{S_1,\ldots,S_n\}$ be a configuration of simplices. Let $F$ be a facet of one of the simplices in the configuration. We say that $F$ is an {\it open} facet if $F$ belongs to only one simplex in the configuration.
\end{defin}

Now we are ready to complete the goal of this section and prove that the Delaunay triangulation of $Q$ must contain only facets of three types.

\begin{lem}
\label{lem:3types}
The Delaunay triangulation of $Q$ does not contain facets of Type D or Type~E.
\end{lem}
\begin{proof} We use the following Depth-First Search algorithm that starts from a pair of adjacent simplices, and add simplices to an open facet of existing configuration while trying to either find frame violations or to certify new configurations by proving that maximal determinant of any positive definite quadratic form with covering radius at most $1$ that contains all simplices of the configuration is below the $4.395397$ threshold.
    
We start from a configuration $\mathcal C_0=\{\sigma_0,\sigma\}$ of two adjacent simplices of Type D or Type E. Our initial stack $\mathfrak S$ is this configuration paired with the quadratic form $Q_0$ (given as $6\times 6$ {\tt float} matrix) that maximizes determinant of positive definite quadratic forms with circumradii at most 1 on simplices of $\mathcal C_0$; this form is obtained numerically using the {\tt maxdet} function by optimizing the concave function $\log \det Q$ over the convex set of all relevant positive definite quadratic forms $Q$.

While the stack is non-empty, we pop the last configuration $\mathcal C_i$ from the stack with the corresponding form $Q_i$, and pick its open facet $F_i$. The facet $F_i$ is chosen among all open facets of the configuration as follows.
\begin{enumerate}
\item For each open facet, we choose all potential apexes of new simplices provided that the open facet will become one of the five Types A, B, C, D, or E.
\item For the potential apexes, we choose only the simplices that satisfy conditions (a) and (b) of Definition~\ref{def:conf}.
\item For all such simplices, we compute approximate circumradii with respect to $Q_i$ and choose the open facet that maximizes the minimum of circumradii. If this maxmin circumradius is greater than $1$, then adding any of these potential simplices will decrease {\tt maxdet} of all potential new configurations.
\item In case there are several open facets with the same maxmin circumradius, we choose the facet with fewer potential simplices.
\end{enumerate}

Once we choose the open facet $F_i$, for each new potential simplex $S_{ij}$ across $F_i$, we check whether adding the class $S_{ij}$ to $\mathcal C_i$ violates the frame condition or whether the new configuration $\mathcal C_i \cup \{S_{ij}\}$ can be certified to have maximal determinant of all positive definite quadratic forms with circumradii of all simplices at most $1$ smaller than $4.395397$. In any of these cases, we know that $\mathcal C_i \cup \{S_{ij}\}$ cannot appear in the Delaunay triangulation of $Q$.

If none of these two conditions hold, we add the new {\it child} configuration $\mathcal C_m:=\mathcal C_i \cup \{S_{ij}\}$ to the stack $\mathfrak S$ together with the numerically computed form $Q_m$ from {\tt maxdet} during the certification attempt. 

Once the stack is empty, we conclude that the initial configuration in the stack cannot appear in the Delaunay triangulation of $Q$ with $\det Q> 4.395397$ and covering radius $1$. If it does appear, then we would reconstruct it in our search after getting a configuration with 360 classes of simplices and no open facets as at every step we pick an open facet, create new configurations by making this open facet of one of the five types, and either certify new configurations or add them to the stack for further exploration, and if we add simplices of an actual triangulation, we cannot certify that it has determinant below $4.395397$.

Overall, the output of the script {\tt 02_step2_types/03_types_DE.py} is the following. Both initial configurations produce empty stacks with the total running time less than one minute\footnote{The time may depend on the environment used. We used a MacBook Pro with a 14-core Apple M4 processor; only single-CPU computations were used.}. 

{\footnotesize\tt\begin{center}
\begin{tabular}{l}
Working on type [2, 1, 1, 0, -1, -1]\\
Total visited: 0 . Stack size:  1 . Working on these roots:  4e-06 seconds\\
Total visited: 10 . Stack size:  3 . Working on these roots:  10.795487 seconds\\
Total visited: 20 . Stack size:  2 . Working on these roots:  18.698048 seconds\\
Total visited for [2, 1, 1, 0, -1, -1] : 23 \\
Total certified: 485 \\
Total children: 649 \\
Total excluded by frame cutoff: 142 \\
Maximal configuration size: 5 \\
Total time spent:  21.082999 \\
\\
Working on type [1, 1, 1, 1, 0, -2] \\
Total visited: 0 . Stack size:  1 . Working on these roots:  4e-06 seconds \\
Total visited for [1, 1, 1, 1, 0, -2] : 7 \\
Total certified: 694 \\
Total children: 881 \\
Total excluded by frame cutoff: 181 \\
Maximal configuration size: 3 \\
Total time spent:  22.569908 \\
\end{tabular}
\end{center}}
\end{proof}

At this point we know that if there is a generic positive definite quadratic form $Q$ with $\Theta(Q)<2.4649$, then the Delaunay triangulation $\mathcal D(Q)$ is unimodular and all facets are of Types A, B, or C.

\section{Facets of Type A}
\label{sec:type-A}

In this section we prove that every simplex in the Delaunay triangulation of the form $Q$ with $\Theta(Q)<2.4649$ has exactly four facets of Type A and that together the simplex and the four neighbors form a very distinctive pattern. Moreover, extending these arguments to further facets of the neighbors, we get that each Type-A connected component in the Delaunay triangulation of $Q$ is formed by subdividing a free sum of three triangles into nine simplices.

These results are obtained by a computational verification of the claims below. The computations are performed in the same way we discarded Type D and Type E facets in Lemma~\ref{lem:3types} with an additional check that traces all configurations we visited. If during the verification of a claim the search encounters a configuration of simplices $\mathcal C'$ which is affinely equivalent to a configuration of simplices $\mathcal C$ that was visited during the same claim, then we mark the configuration $\mathcal C'$ as resolved. Since our search only adds simplices to every configuration, searching through one affine copy is sufficient for both.

Our affine symmetry check comes in two different options depending on the type of the claim we are proving; for each claim, we specify which symmetry check we are using. Suppose we are given two simplex configurations $\mathcal C$ and $\mathcal C'$.
\begin{itemize}
    \item In the {\it unrestricted} case, we try to find any affine map that maps one configuration to the other. To do that, we pick a class of simplices $S$ in $\mathcal C$ and map it to the class of the standard simplex $\sigma_0$ by an affine map. Among all such maps we choose the one that gives the lexicographically smallest image of selected simplices from each class. For the image of each class, we choose the simplex where the origin is lexicographically the smallest vertex and pick the smallest of the two centrally symmetric copies in each class. We repeat that for $\mathcal C'$ and conclude that two configurations are affinely equivalent if and only if they give the same output.

    \item In the {\it restricted} case, the process is similar, but we fix one simplex class in $\mathcal C$ and one simplex class in $\mathcal C'$ in advance and these classes must be mapped to each other (or must be mapped to $\sigma_0$ in our lexicographic search). 
\end{itemize}

The rest of the section contains the claims we prove and some details about the computations; the complete logs for each computation are contained in the supplied package \cite{package}. We start by proving that a single Type A facet implies a fixed structure of nine simplices; this is done in Claims~\ref{claim1} and \ref{claim2} and Corollary~\ref{cor:free-sum}.

As before, we assume that $Q$ is a generic positive definite quadratic form with $\det Q > 4.395397$ and covering radius $1$. In the previous sections we already established that all facets of the Delaunay triangulation $\mathcal D(Q)$ are of Types A, B, or C.

\begin{claim}
\label{claim1}
Suppose the Delaunay triangulation of $Q$ contains the standard simplex $\sigma_0$ and its neighbor
$$\sigma_1 = \conv \left(\begin{array}{l}
(1,1,1,-1,0,0)^t,\\
(1,0,0,0,0,0)^t,\\
(0,1,0,0,0,0)^t,\\
(0,0,1,0,0,0)^t,\\
(0,0,0,1,0,0)^t,\\
(0,0,0,0,1,0)^t,\\
(0,0,0,0,0,1)^t
\end{array}\right)$$
across the common facet which is a Type A facet in this setting.
Then it also contains their common neighbor
$$\sigma_2 = \conv \left(\begin{array}{l}
(0,0,0,0,0,0)^t,\\
(1,0,0,0,0,0)^t,\\
(0,1,0,0,0,0)^t,\\
(0,0,1,0,0,0)^t,\\
(1,1,1,-1,0,0)^t,\\
(0,0,0,0,1,0)^t,\\
(0,0,0,0,0,1)^t
\end{array}\right)$$
and the common facets of $\sigma_2$ with $\sigma_0$ and $\sigma_1$ are both of Type A.
\end{claim}

\begin{proof}
The proof is obtained by running {\tt python} script {\tt 03_step3_typeA/claim1.py} from the package \cite{package}; the output is contained in {\tt 03_step3_typeA/claim1.log}.

The idea is as follows. Once we fix the two simplices $\sigma_0$ and $\sigma_1$, we consider the list $\mathcal L$ of potential neighbors of $\sigma_0$ in the Delaunay triangulation of $Q$ across the facet opposite to $(0,0,0,1,0,0)$. This list is finite because only facets of Types A, B, or C are possible. Let $\mathcal L'=\mathcal L\setminus \{\sigma_2\}$ be the list of potential neighbors across the facet of $\sigma_0$ with $\sigma_2$ removed.

We start from the stack $\mathfrak{S}:=\{\{\sigma_0,\sigma_1,S\}: S\in \mathcal L'\}$ of simplex configurations and run our search adding potential simplices to open facets selecting open facets in the same way as in Lemma~\ref{lem:3types}. Once every configuration from the stack is removed, we conclude that $\sigma_2$ must be in the Delaunay triangulation of $Q$ provided $Q$ satisfies our usual restrictions on the determinant and covering radius because for each other case we arrive to an upper bound on $\det Q$ below $4.395397$.

Overall, the script identifies 23 potential neighbors other than $\sigma_2$ (after checking the frame conditions and certifying neighbors removed by {\tt maxdet}) in 7 unrestricted symmetry classes. The script discards all 7 classes in about 1 hour and 49 minutes with the largest configuration reaching size 13 and largest case visiting 5916 configurations (for which we needed to compute potential neighbors). 

We note that the configuration with potential neighbor $S$ with vertex $(1,1,0,-1,0,0)^t$ and initial {\tt maxdet} $8.643850173452133$ of the configuration of three simplices $\{\sigma_0,\sigma_1,S\}$ took about 1 hour and 33 minutes of the computation (almost 85\% of the total time) while the maximal configuration size was reached for another potential neighbor.
\end{proof}

Once we find these three simplices, we can find more neighbors of Type A of the standard simplex $\sigma_0$.

\begin{claim}
\label{claim2}
In the setting of Claim~\ref{claim1}, the simplex 
$$\sigma_3 = \conv \left(\begin{array}{l}
(0,0,0,0,0,0)^t,\\
(1,0,0,0,0,0)^t,\\
(0,1,0,0,0,0)^t,\\
(0,0,1,0,0,0)^t,\\
(0,0,0,1,0,0)^t,\\
(1,1,1,0,-1,-1)^t,\\
(0,0,0,0,0,1)^t
\end{array}\right)$$ is also present in the Delaunay triangulation of $Q$. This simplex is another neighbor of $\sigma_0$ and the common facet of $\sigma_0$ and $\sigma_3$ is of Type A.
\end{claim}

\begin{proof}
The proof is obtained by running {\tt python} script {\tt 03_step3_typeA/claim2.py} from the package \cite{package}; the output is contained in {\tt 03_step3_typeA/claim2.log}.

The idea is similar to the previous claim. The initial configuration consists of three classes of simplices $\{\sigma_0,\sigma_1,\sigma_2\}$. We consider all potential neighbors of $\sigma_0$ across the facet opposite to $(0,0,0,0,1,0)^t$ and the stack contains all configurations with the neighbors different from our target simplex $\sigma_3$. We discard all potential configurations either through the frame condition or through {\tt maxdet} certificate.

The script identifies 30 potential neighbors other than $\sigma_3$ in 10 unrestricted symmetry classes. The script discards all 10 classes in about 1 hour and 12 minutes with the largest configuration reaching size 16 and the largest case visiting 1574 configurations.
\end{proof}

\begin{cor}
\label{cor:free-sum}
In the setting of Claim~\ref{claim1}, the Delaunay triangulation of $Q$ contains a triangulation of the free sum of three triangles $T_1,T_2,T_3$ with the following vertices.
\begin{center}
\begin{tabular}{rlll}
$T_1$:& $(1,0,0,0,0,0)^t,$ & $(0,1,0,0,0,0)^t,$ & $(0,0,1,0,0,0)^t,$ \\
$T_2$:& $(0,0,0,0,0,0)^t,$ & $(0,0,0,1,0,0)^t,$ & $(1,1,1,-1,0,0)^t,$ \\
$T_3$:& $(0,0,0,0,1,0)^t,$ & $(0,0,0,0,0,1)^t,$ & $(1,1,1,0,-1,-1)^t$. \\
\end{tabular}
\end{center}
The free sum is triangulated into nine simplices with each simplex containing all three vertices of $T_1$ and two vertices of $T_2$ and $T_3$ each.
\end{cor}

\begin{proof}
We note that Claims~\ref{claim1} and \ref{claim2} already give four of the required nine simplices. Applying Claim~\ref{claim1} to simplices $\sigma_0$ and $\sigma_3$, we get that the Delaunay triangulation of $Q$ contains simplex $\sigma_4$, the convex hull of the nine points except $(0,0,0,0,0,1)^t$ and $(1,1,1,-1,0,0)^t$. So the simplex $\sigma_0$ contains four Type A facets with neighbors within this free sum.

The remaining four simplices can be obtained by applying Claims~\ref{claim2} and~\ref{claim1} to the other Type A facets within this free sum of three triangles. At this point $\sigma_0$ has common facets with $\sigma_1,\sigma_2,\sigma_3$, and $\sigma_4$, while $\sigma_1$ has common facets only $\sigma_0$ and $\sigma_2$ from these five simplices. Applying an affine transformation that swaps $\sigma_1$ and $\sigma_0$ and preserves the lattice, we can find two more Type A facets of $\sigma_1$ using Claims~\ref{claim2} and~\ref{claim1}. Similarly, we can find two more Type A facets for the other simplices in this cluster.

Overall, we get nine simplices forming the claimed triangulation of the free sum and each of the nine simplices has exactly four Type A facets with neighbors within the free sum.
\end{proof}

Particularly, Corollary~\ref{cor:free-sum} shows that if a simplex has one Type A facet, then it has at least four Type A facets and the four neighbors are defined by the first neighbor. Our next claim shows that none of the remaining three facets of the simplex can be of Type A.

\begin{claim}
\label{claim3}
In the setting of Claim~\ref{claim1} and Corollary~\ref{cor:free-sum}, the three facets of $\sigma_0$ opposite to vertices $(1,0,0,0,0,0)^t$, $(0,1,0,0,0,0)^t$, and $(0,0,1,0,0,0)^t$ are all not of Type A.
\end{claim}

\begin{proof}
Due to the symmetry of the whole configuration of nine simplices from Corollary~\ref{cor:free-sum} under permutations of the first three coordinates, it is enough to show that the facet of $\sigma_0$ opposite to $(1,0,0,0,0,0)^t$ is not of Type A.

The proof is obtained by running {\tt python} script {\tt 03_step3_typeA/claim3.py} from the package \cite{package}; the output is contained in {\tt 03_step3_typeA/claim3.log}.

The idea is similar to the previous claims but now we consider only potential neighbors assuming the facet is of Type A. The initial configuration consists of five classes of simplices $\{\sigma_0,\sigma_1,\sigma_2,\sigma_3,\sigma_4\}$\footnote{We could have started from the nine simplices from Corollary~\ref{cor:free-sum}, but the computation is already very small so we avoided running {\tt maxdet} on larger configurations.}. We consider all potential neighbors of $\sigma_0$ across the facet opposite to $(1,0,0,0,0,0)^t$ provided this facet is of Type A and discard all potential configurations either through the frame condition or through {\tt maxdet} certificate.

The script identifies 8 potential neighbors in one unrestricted symmetry class. The script discards the class in less than five seconds by finding an open facet where adding any of the 7 potential simplices will give configuration with certified maximal determinant below $4.395397$.
\end{proof}

Thus the first three claims establish that once a unimodular simplex $\sigma$ in the Delaunay triangulation $\mathcal D(Q)$ has one Type A facet, then it belongs to a triangulated free sum of three triangles (with each simplex in the triangulation having four Type A facets with neighbors within the free sum), and all neighbors of $\sigma$ outside the free sum are of Types B or C.

Our next goal is to show that every simplex in the Delaunay triangulation of $Q$ has a Type A facet. This will guarantee that the Delaunay triangulation of $Q$ can be split into a collection of disjoint free sums of triangles that are further triangulated. The existence of Type A facet for each simplex in $\mathcal D(Q)$ is established in the next claim.

\begin{claim}
\label{claim4}
Each simplex in the Delaunay triangulation of $Q$ has a Type A facet.
\end{claim}

\begin{proof}
Suppose that the Delaunay triangulation of $Q$ contains a simplex without Type A facets. We may assume that this simplex is the standard simplex $\sigma_0$.

The rest of the proof is obtained by running {\tt python} script {\tt 03_step3_typeA/claim4.py} from the package \cite{package}; the output is contained in {\tt 03_step3_typeA/claim4.log}.

We start from the configuration of the standard simplex $\{\sigma_0\}$. For one facet, we consider all potential neighbors across that facet assuming the facet is of Type B or Type C. We create new configurations by adding new simplices that do not violate the frame condition and do not give a configuration with maximal determinant of relevant positive definite quadratic forms $Q$ below $4.395397$. This process is repeated sequentially for all seven facets of $\sigma_0$. We also discard configurations for which we already added one equivalent under symmetries preserving $\sigma_0$ as we need to surround this particular simplex by neighbors across facets of Types B or C only.

At the end we get a stack of 415 non-equivalent configurations of eight simplices each ($\sigma_0$ and its seven neighbors) with maximum {\tt maxdet} $7.192904000667755$ among these configurations. The script then allows all three types of facets for all other open facets within each configuration in the stack and keeps a common collection of all visited configurations that are compared up to symmetries without restrictions.

The script depletes the stack of all 415 cases in about 1 hour and 50 minutes (including the initial setup of configurations of $\sigma_0$ and its neighbors) with the largest configuration reaching size 14 and visiting 4042 configurations for all cases combined.
\end{proof}

Altogether, Claims~\ref{claim1},~\ref{claim2},~\ref{claim3},~\ref{claim4}, and Corollary~\ref{cor:free-sum} prove that the simplices of the Delaunay triangulation of $Q$ can be combined into groups of nine such that
\begin{itemize}
\item The nine simplices in every group form a free sum of three triangles;
\item Each simplex has exactly four Type A facets with all neighbors in its group; and
\item Each two neighbors from different groups share facet of Type B or Type C.
\end{itemize}
These properties are shared by the Delaunay triangulation of $Q_6^{c}$ and by any triangulation that subdivides the Delaunay tessellation of $\mathsf E_6^*$ as shown in Lemma~\ref{lem:typesABC}, and our final step is to show that the Delaunay triangulation of $Q$ is such a subdivision.

\section{The final step}
\label{sec:e6star}

As the final step we prove that the Delaunay triangulation of $Q$ subdivides the Delaunay tessellation of $\mathsf E_6^*$. This allows us to use Theorem~\ref{thm:SV} by Sch\"urmann and Vallentin and complete the proof of our main Theorem~\ref{thm:main}. At this point, we may assume that $Q$ is a generic positive definite quadratic form with $\det Q > 4.395397$ and covering radius at most $1$, and we use all structural properties of Type A, Type B, and Type C facets from the previous section.

\begin{claim}
\label{claim-adj-sums}
Suppose that the Delaunay triangulation of $Q$ contains nine unimodular simplices subdividing the free sum $P_{123}$ of three triangles $T_1$, $T_2$, and $T_3$ with vertices given by
\begin{center}
\begin{tabular}{rlll}
$T_1$:& $(1,0,0,0,0,0)^t,$ & $(0,1,0,0,0,0)^t,$ & $(0,0,1,0,0,0)^t,$ \\
$T_2$:& $(0,0,0,0,0,0)^t,$ & $(0,0,0,1,0,0)^t,$ & $(1,1,1,-1,0,0)^t,$ \\
$T_3$:& $(0,0,0,0,1,0)^t,$ & $(0,0,0,0,0,1)^t,$ & $(1,1,1,0,-1,-1)^t$. \\
\end{tabular}
\end{center}
Then this Delaunay triangulation contains nine unimodular simplices subdividing the free sum $P_{456}$ of three triangles $T_4$, $T_5$, and $T_6$ with vertices
\begin{center}
\begin{tabular}{rlll}
$T_4$:& $(-1,-1,-1,1,1,1)^t,$ & $(0,1,0,0,0,0)^t,$ & $(0,0,1,0,0,0)^t,$ \\
$T_5$:& $(0,0,0,0,0,0)^t,$ & $(0,0,0,1,0,0)^t,$ & $(-1,0,0,0,1,1)^t,$ \\
$T_6$:& $(0,0,0,0,1,0)^t,$ & $(0,0,0,0,0,1)^t,$ & $(-1,0,0,1,0,0)^t$. \\
\end{tabular}
\end{center}
\end{claim}

\begin{proof}
We first notice that all three triangulations of the free sum $P_{123}=\conv\{T_1,T_2,T_3\}$ can be mapped to each other while preserving the lattice $\Z^6$ and the free sum $P_{456}=\conv \{T_4,T_5,T_6\}$. Therefore, we can fix one triangulation of the first free sum and assume that these nine simplices are in $\mathcal D(Q)$; we will assume that nine simplices $\{\sigma_0,\ldots,\sigma_8\}$ from Corollary~\ref{cor:free-sum}, including the standard simplex $\sigma_0$ are in $\mathcal D(Q)$.
    
The rest of the proof is split into two computations both realized in the script {\tt 04_step4_cells/free_sum_cells.py} in the package \cite{package}. The output is in {\tt 04_step4_cells/free_sum_cells.log}.

First, we look at the facet $F=\conv \{\mathbf 0, e_2,e_3,e_4,e_5,e_6\}$ of $\sigma_0$ opposite to $(1,0,0,0,0,0)^t$. The neighbor of $\sigma_0$ across the facet $F$ is outside of $P_{123}$, so $F$ has to be of Type B or Type C. We test all potential simplices $\sigma$ that can be this neighbor by checking the frame conditions and {\tt maxdet} for configurations $\{\sigma_0,\ldots,\sigma_8,\sigma\}$. Among all the candidates, only three potential simplices pass the threshold allowing $\det Q> 4.395397$ and these simplices are exactly simplices within $P_{456}$:
$$\conv\{F, (-1,-1,-1,1,1,1)^t\}, \quad \conv\{F, (-1,0,0,0,1,1)^t\}, \quad  \conv\{F, (-1,0,0,1,0,0)^t\}.$$

After that, for each of the three remaining candidates $\sigma$, our goal is to show that all nine simplices of the triangulation of $P_{456}$ containing $\sigma$ are in $\mathcal D(Q)$; if we use symmetries of the configurations of 10 simplices (the triangulation of $P_{123}$ and $\sigma$), then it is enough to consider only two cases. To treat each case, we choose one facet $F'$ of $\sigma$ inside $P_{456}$ and another simplex $\sigma'$ with facet $F'$ in $P_{456}$. Our script checks that among all possible neighbors of $\sigma$ across $F'$, only $\sigma'$ allows determinant of $Q$ to stay above $4.395397$. This is done similarly to our proof of Claim~\ref{claim1}. 

Note that we allow $F'$ to be of any of the three Types A, B, or C, but our script confirms that only a specific adjacent simplex makes determinant above $4.395397$ and that simplex turns $F'$ into Type A facet. More precisely, 
\begin{itemize}
\item If $\sigma = \conv\{F, (-1,-1,-1,1,1,1)^t\}$, then we show that simplex 
$$\sigma'=\conv \left(\begin{array}{l}
(0,0,0,0,0,0)^t,\\
(-1,-1,-1,1,1,1)^t,\\
(0,1,0,0,0,0)^t,\\
(0,0,1,0,0,0)^t,\\
(0,0,0,1,0,0)^t,\\
(0,0,0,0,1,0)^t,\\
(-1,0,0,1,0,0)^t
\end{array}\right)$$
must be present in the Delaunay triangulation of $Q$,
\item If $\sigma = \conv\{F, (-1, 0, 0, 0, 1, 1)^t\}$, then we show that simplex 
$$\sigma'=\conv \left(\begin{array}{l}
(0,0,0,0,0,0)^t,\\
(-1, 0, 0, 0, 1, 1)^t,\\
(0,1,0,0,0,0)^t,\\
(0,0,1,0,0,0)^t,\\
(0,0,0,1,0,0)^t,\\
(0,0,0,0,1,0)^t,\\
(-1,0,0,1,0,0)^t
\end{array}\right)$$
must be present in the Delaunay triangulation of $Q$, and 
\item If $\sigma = \conv\{F, (-1,0,0,1,0,0)^t\}$, then there is a symmetry of that case to the previous one preserving $P_{123}$ and $P_{456}$; in both cases the facet $F$ is of Type B and this symmetry is recognized by script's check of restricted symmetries fixing $\sigma_0$.
\end{itemize}

After that, we can apply Claims~\ref{claim1} and~\ref{claim2}, and Corollary~\ref{cor:free-sum} to simplices $\sigma$ and $\sigma'$ as their common facet $F'$ is of Type A, and we get that the Delaunay triangulation of $Q$ contains nine simplices triangulating $P_{456}$.

Overall, the two-step script that searches potential neighbors $\sigma$ of $\sigma_0$ and discards all neighbors of candidates $\sigma$ across $F'$ takes less than one minute and visits 15 different configurations with maximal size 15; note that the triangulation of $P_{123}$ with added $\sigma$ and its neighbor $\sigma'$ already give $11$ simplices.
\end{proof}

The claim above allows us to apply Theorem~\ref{thm:SV} as described in the proof of Theorem~\ref{thm:main} in Section~\ref{sec:plan}. Indeed, it guarantees that once we have one triangulated free sum of three triangles in $\mathcal D(Q)$, its neighboring free sum clusters in $\mathcal D(Q)$ are defined uniquely (but triangulation of each specific free sum is not uniquely defined). We have shown this property for one specific adjacent free sum in Claim~\ref{claim-adj-sums}, and the rest follow from symmetries of the lattice $\Z^6$ preserving the free sum $P_{123}$ which acts transitively on its facets.

Therefore, one free sum gives us all free sums provided $\det Q>4.395397$ and $R(Q)\leq 1$. Since the properly scaled form $Q_6^c$ satisfies these properties, the free sums for $Q$ and for $Q_6^c$ coincide after an appropriate transformation from $GL_6(\Z)$. Thus, the Delaunay triangulation of $Q$ subdivides the Delaunay tessellation of $\mathsf E_6^*$ and Theorem~\ref{thm:SV} and Proposition~\ref{prop:SV-unique} by Sch\"urmann and Vallentin guarantee that if $\Theta(Q)\leq \Theta(Q_6^c)$, then $Q$ is in the secondary cone of $Q_6^c$ or in the equivalent cone. The section below provides the remaining piece for the proof of Theorem~\ref{thm:main} for generic forms by giving an exact expression for $Q_6^c$ and proving that it is unique minimizer for the covering density in the corresponding secondary cone.

\section{Exact expression for \texorpdfstring{$Q_6^c$}{Q6c} and its local optimality}
\label{sec:exact}

In this section we give details about the exact quadratic form $Q_6^c$ and provide details on how to check that this quadratic form has covering density below $2.4649$. We also prove that this form is the unique minimizer for covering density within closure of its secondary cone. Our verification can be done by running {\tt SageMath} script {\tt exact_form/Q6c_exact.sage} in the package \cite{package}. Here are the details of the verification process and we prove the minimizing property after that.

Recall from Section~\ref{sec:lattice} that
\begin{equation}\label{eq:exact}
Q_6^c = Q_{\mathsf E_6^*} + x R_2 + y R_3 + z R_4
\end{equation} where approximate values for the three parameters found by Sch\"urmann and Vallentin in \cite{SV06} are $x = 0.15266328480099$, $y =
0.32884740614948$, and $z = 0.13827491241153$ and the matrices for the four quadratic forms in the right-hand side are given in Section~\ref{sec:lattice}. The precise values for the parameters are defined by the following polynomials in $\Z[t]$.

The minimal polynomial for $x$ is 
{\tiny
\begin{center}    
\begin{tabular}{l}
$f_x = 2402475394684564224t^{19} - 6437459726723417216t^{18} - 19844323601779568936t^{17} -17859560598939121702t^{16} + $ \\ $+14250122198850976100t^{15} + 41165222327347137567t^{14} + 49216175347976933090t^{13} + 43158308490593567534t^{12} +$ \\ $+26232029904300449169t^{11} + 8362849663056498774t^{10} - 818523219908212395t^9 - 2195528505965172165t^8 -$\\$- 1031183161797013680t^7 - 237999981714295380t^6 - 21679790140751280t^5 + 2352381933354240t^4 +$\\$+ 864877778092800t^3 + 95869730246400t^2 + 4714527974400t + 85764096000 $
\end{tabular}
\end{center}}
and $x$ is the only real root of that polynomial on the interval $(0.152,0.153)$. The minimal polynomial for $y$ is 
{\tiny
\begin{center}
\begin{tabular}{l}
$f_y = 25025785361297544t^{19} + 110659621549017256t^{18} + 207579230411121962t^{17} + 203875308944320384t^{16} + $\\ $+ 93088201565148501t^{15} - 13033553733693743t^{14} - 46680879865229733t^{13} - 35612559350035656t^{12} - $ \\ $-20151647881070942t^{11} - 10240200149614773t^{10} - 3985625973319548t^9 - 778260069690204t^8 + $\\$+142625553848517t^7 + 150790962060339t^6 + 50469229897071t^5 + 9714710419044t^4 +$\\$+ 1147603507980t^3 + 80381041065t^2 + 3040751160t + 54928800$
\end{tabular}
\end{center}
}
and $y$ is the only root of this polynomial on the interval $(0.328,0.329)$. The minimal polynomial for $z$ is 
{\tiny
\begin{center}
\begin{tabular}{l}
$f_z = 33367713815063392t^{19} + 298184979369241392t^{18} + 1018678597166385148t^{17} + 1792863305110275168t^{16} +$\\$+ 1692369572666137417t^{15} + 682403306467785788t^{14} - 170679511359803597t^{13} - 272603474894126335t^{12} -$\\$- 39113760696426042t^{11} + 62641854799250980t^{10} + 26612012992402650t^9 - 10764655908844220t^8 -$\\$- 15715746862076160t^7 - 8495819933393120t^6 - 2805313126456480t^5 - 573430106753600t^4 -$\\$- 46647810304000t^3 + 8656969766400t^2 + 2143805184000t + 88197120000$
\end{tabular}
\end{center}
}
and $z$ is the only root of this polynomial on the interval $(0.138,0.139)$. The three polynomials above were found using {\tt ChatGPT} and {\tt Claude}.

The following are slightly better approximations for $x$, $y$, and $z$ compared to \cite{SV06}.
\begin{center}
\begin{tabular}{l}
$x\approx 0.1526632848049188$\\
$y\approx 0.3288474061128059$\\
$z\approx 0.1382749123683778$
\end{tabular}
\end{center}

We note that both $y$ and $z$ are in the field $\Q(x)$ but the corresponding rational expressions involve numerators and denominators with more than 200 digits, so we are not presenting them here; the relevant check is done in the accompanying package. Moreover, it is much faster to perform computations when all values are taken from $\Q(x)$ below. The description in $\Q(x)$ we discuss below was found by {\tt Claude}. The computational verification was independently implemented in {\tt SageMath} by the author.

For our verification, we check that the following {\tt Claude} supplied data satisfies the following properties.
\begin{enumerate}
\item We start from 40 classes of free sums and check that they form the known Delaunay tessellation of the quadratic form $Q_{\mathsf E_6^*}$.
\item For each free sum, the {\tt Claude} data contains description of its triangulation, and we verify that this is the Delaunay triangulation of the free sum with respect to the quadratic form $Q_6^c$ from \eqref{eq:exact}. This is done by checking local Delaunay property for each pair of adjacent simplices within the free sum. We don't need to check other adjacent simplices as the $\mathcal D(Q_6^c)$ as given in \eqref{eq:exact} subdivides $\mathcal D(\mathsf E_6^*)$ by Lemma~\ref{lem:closure}.
\item We check that all simplices in the Delaunay triangulation have squared circumradii that do not exceed a pre-computed squared circumradius $\mu \in \Q(x)$ and $\mu$ is actually the squared covering radius for the quadratic form $Q_6^c$.
\item For the list of 180 simplices provided by {\tt Claude}, we check that each simplex is in $\mathcal D(Q_6^c)$ and has the confirmed maximal circumradius.
\item As the last step, we run the certification process from Section~\ref{sec:lemmas}. The 180 simplices with the maximal circumradius are split into three orbits under the symmetry group of $Q_6^c$, and each orbit has its own rank 1 $7\times 7$ matrix with entries from $\Q(x)$; the order of vertices is important for that as reordering changes the matrix $B_L$ from Lemma~\ref{lem:B_L}. We verify that these three matrices are positive semidefinite, and that 
\begin{itemize}
    \item the matrix $P_\mathcal L$ is positive definite, and
    \item the right-hand side $C_\mathcal L$ is positive,
\end{itemize}
where $\mathcal L$ is the collection of 180 simplices. The resulting certified upper bound on determinant of positive definite form $Q$ coincides with determinant of properly scaled $Q_6^c$. This implies that $Q_6^c$ given by \eqref{eq:exact} minimizes the covering density in its secondary cone as claimed.
\end{enumerate}

Overall, the verification process takes less than 15 minutes and majority of the time is spent on computing circumradii of simplices. All computations are exact and are performed in $\Q(x)$.

We also compute the precise covering density as
$$\Theta(Q_6^c)=\frac{R^6(Q_6^c)}{\sqrt{\det Q_6^c}}\kappa_6 \approx 2.46480067066519$$
as squared covering radius and determinant are computed as elements of $\Q(x)$ during the {\tt SageMath} verification process and $\kappa_6=\frac{\pi^3}6$.

Our final step is to show that $Q_6^c$ is the only minimizer of the covering density in the closure of its secondary cone.

\begin{lem}
\label{lem:exact}    
The quadratic form $Q_6^c$ from \eqref{eq:exact} is the unique (up to scaling) minimizer of the covering density in the closure of its secondary cone.
\end{lem}

\begin{proof}
Let $Q$ be another minimizer. We can scale it so the covering radius of $Q$ is $1$. Since Delaunay triangulation of $Q_6^c$ subdivides Delaunay tessellation of $Q$, all 180 simplices used for the certification process above will have circumradii at most $1$ with respect to $Q$ as well. If we follow the same certification process, then we see that the AM-GM inequality 
$$\det Q \cdot \det P_{\mathcal L} = \det(AP_{\mathcal L}A)\leq \left(\frac{C_{\mathcal L}}{6}\right)^6$$
where $A=Q^{\frac12}$ from Section~\ref{sec:lemmas} must turn into an equality as it did for $Q_6^c$. This means that all eigenvalues of $A P_{\mathcal L}A$ must be equal and $Q$ is a scaled matrix $P_{\mathcal L}^{-1}$. Since $P_{\mathcal L}$ depends only on the collection of simplices and the matrices $N_L$, $Q$ is a scaled form $Q_6^c$ as claimed.
\end{proof}

\section{Experimental results in \texorpdfstring{$\R^7$}{R7}}

One of the most reasonable questions is whether a similar approach can be carried out in other dimensions. It looks straightforward for $d=2,3,4$, and $5$ where the best lattice sphere coverings are given by the lattices $\mathsf A_d^*$ and all facets are analogs of our Type B. Moreover, for $d=5$, we need to consider up to 60 classes of simplices and though our tools are tailored specifically to the six-dimensional case, they can handle configurations of that size. So we expect that similar approaches to local combinatorics are possible for the known cases $d\leq 5$.

For higher dimensions, and specifically for $d=7$, the currently best known lattice with the smallest covering density is the lattice $\mathsf L_7^c$ found by Sch\"urmann and Vallentin in \cite{SV06}. This lattice has Gram matrix 
$$
Q_7^c = \left(
\begin{array}{ccccccc}
12& 1& 1& 1& 1& 1& 5\\ 
1& 12& 1& 1& 1& 1& 5\\ 
1& 1& 12& 1& 1& 1& 5\\ 
1& 1& 1& 12& 1& 1& 5\\ 
1& 1& 1& 1& 12& 1& 5\\ 
1& 1& 1& 1& 1& 12& -6\\ 
5& 5& 5& 5& 5& -6& 14
\end{array}
\right)
$$
and covering density $\Theta(\mathsf L_7^c)\approx 2.900024$. 

We asked {\tt Claude} to try a similar approach for this lattice. First of all, the Delaunay triangulation $\mathcal D(Q_7^c)$ also contains three types of facets, seven-dimensional counterparts of Types A, B, and C. 
\begin{center}
\begin{tabular}{rl}
Type $A_7$: & (1,1,1,0,0,0,-1)\\
Type $B_7$: & (1,1,0,0,0,0,0)\\
Type $C_7$: & (1,1,1,1,0,-1,-1)
\end{tabular}
\end{center}

But the process of exclusion for all additional types seems to require much more computations or an additional idea. For example, {\tt Claude} identified a configuration of more than 600 simplices with six facet types that are not refuted by the {\tt maxdet} function.


\begin{thebibliography}{20}

\bibitem{AS26} S.~Ahrend, A.~Sch\"urmann, Computational Approaches to Lattice Packing and Covering Problems via C-types, {\it in preparation}.

\bibitem{Bam54} R.~P.~Bambah, On lattice coverings by spheres, {\it Proc. Nat. Inst. Sci. India}, {\bf 20} (1954), 25--52.

\bibitem{BS82}  R.~P.~Bambah and N.~J.~A.~Sloane, On a problem of Ryshkov concerning lattice coverings, {\it Acta Arith.}, {\bf 42}:1 (1982), 107--109, \url{http://eudml.org/doc/205857}.


\bibitem{Ble62} M. N. Bleicher, Lattice coverings of $n$-space by spheres, {\it Canad. J. Math.}, 14 (1962), 632--650, \url{https://doi.org/10.4153/CJM-1962-054-7}.

\bibitem{Bli35} H. F. Blichfeldt, The minimum values of positive quadratic forms in six, seven and eight variables, {\it Math. Z.}, {\bf 39} (1935), 1--15, \url{https://doi.org/10.1007/BF01201341}

\bibitem{BGLPS25} B.~Bukh, J.~Gao, X.~Liu, O.~Pikhurko, and S.~Sun. Covering large-dimensional Euclidean spaces by random translates of a given convex body, {\it Disc. \& Comp. Geom.}, {\bf 76} (2026), 2424--2442, \url{https://doi.org/10.1007/s00454-026-00871-5}.

\bibitem{CK09} H.~Cohn, A.~Kumar, Optimality and uniqueness of the Leech lattice among lattices, {\it Ann. of Math.}, {\bf 170} (2009), 1003--1050, \url{https://doi.org/10.4007/annals.2009.170.1003}.

\bibitem{CKMRV17} H.~Cohn, A.~Kumar, S.~D.~Miller, D.~Radchenko, and M.~Viazovska. The sphere packing problem in dimension 24, {\it Ann. of Math.}, {\bf 185}:3 (2017), 1017--1033, \url{https://doi.org/10.4007/annals.2017.185.3.8}.

\bibitem{CS88} J.~H.~Conway and N.~J.~A.~Sloane, Sphere packings, lattices and groups, Springer-Verlag, New York, 1988, \url{https://doi.org/10.1007/978-1-4757-6568-7}.

\bibitem{CS91} J.~H.~Conway and N.~J.~A.~Sloane, The Cell Structures of Certain Lattices, In {\it Miscellanea Mathematica}, 1991, 71--107, \url{https://doi.org/10.1007/978-3-642-76709-8_5}.

\bibitem{Del38} B.~N.~Delone, The geometry of positive quadratic forms, {\it Uspekhi Mat. Nauk}, Part I: {\bf 3} (1937), 16--62, \url{https://www.mathnet.ru/eng/rm8903}. Part II: {\bf 4} (1938), 102--164, \url{https://www.mathnet.ru/eng/rm7118}.

\bibitem{DDRS70} B.~N.~Delone, N.~P.~Dolbilin, S.~S.~Ryshkov, and M.~I.~Stogrin, A new construction of the theory of lattice coverings of an $n$-dimensional space by congruent balls, {\it Math. USSR-Izv.}, {\bf 4}:2 (1970), 293--302, \url{https://doi.org/10.1070/IM1970v004n02ABEH000904}.

\bibitem{DR63} B.~N.~Delone and S.~S.~Ryshkov, Solution of the problem on the least dense lattice covering of a 4-dimensional space by equal spheres, {\it Dokl. Akad. Nauk SSSR}, {\bf 152}:3 (1963), 523--524, \url{https://www.mathnet.ru/eng/dan28581}.

\bibitem{DGSW16} M.~Dutour~Sikiri\'c, A.~Garber, A.~Sch\"urmann, and C.~Waldmann, The complete classification of five-dimensional Dirichlet-Voronoi polyhedra of translational lattices, {\it Acta Cryst. A}, {\bf 72}:6 (2016), 673--683, \url{https://doi.org/10.1107/S2053273316011682}.

\bibitem{DR14} M. Dutour Sikiri\'c, K.~Rybnikov, 
Delaunay polytopes derived from the Leech lattice, {\it Journal de Th\'eorie des Nombres de Bordeaux}, {\bf 26}:1 (2014), 85--101, \url{http://eudml.org/doc/275734}.

\bibitem{DSV08} M. Dutour Sikiri\'c, A. Sch\"urmann, and F. Vallentin, A generalization of Voronoi's reduction theory and its application, {\it Duke Math. J.}, {\bf 142}:1 (2008), 127--164, \url{https://doi.org/10.1215/00127094-2008-003}.

\bibitem{DvW25} M. Dutour Sikiri\'c and Wessel van Woerden, The lattice packing problem in dimension 9 by Voronoi's algorithm, {\it arXiv preprint}, arXiv:2508.20719, 2025, \url{https://doi.org/10.48550/arXiv.2508.20719}.


\bibitem{Fed85} E.~S.~Fedorov, Elements of the Study of Figures, Zap. Mineralog. Obsch., 1885.

\bibitem{Gam62}  A. F. Gameckii, On the theory of covering Euclidean $n$-space by equal spheres, {\it Soviet Math. Dokl.}, {\bf 146}:5 (1962), 991--994, \url{https://www.mathnet.ru/eng/dan27080}. 

\bibitem{Gam63} A. F. Gameckii, The optimality of Voronoi's lattice of first type among lattices of first type of arbitrary dimension, {\it Soviet Math. Dokl.}, {\bf 151}:3 (1963), 482--484, \url{https://www.mathnet.ru/eng/dan28324}



\bibitem{GLPS26} J.~Gao, X.~Liu, O.~Pikhurko, and S.~Sun, New upper bound for lattice covering by spheres, {\it Mathematika}, {\bf 72}:e70066 (2026), \url{https://doi.org/10.1112/mtk.70066}.


\bibitem{Gar25} A.~Garber, Voronoi conjecture for five-dimensional parallelohedra, {\it Invent. Math.}, {\bf 240}:2 (2025), 587--631, \url{
https://doi.org/10.1007/s00222-025-01325-0}.

\bibitem{package} A.~Garber, Accompanying {\tt python} and {\tt SageMath} package, available at \url{https://github.com/alexeygarber/6dim-lattice-covering}, version of November 28, 2026.



\bibitem{Gau31} C. F. Gauss, Untersuchungen \"uber die Eigenschaften der positiven tern\"aren quadratischen Formen von Ludwig August Seeber, G\"ottingische gelehrte Anzeigen, July 9, 1831, Reprinted in Werke, Vol. 2, K\"onigliche Gesellschaft der Wissenschaften, G\"ottingen, 1863, 188--196.

\bibitem{Hal05} T.~C.~Hales, A proof of the Kepler conjecture, {\it Ann. of Math.}, {\bf 162}:3 (2005), 1065--1185, \url{https://doi.org/10.4007/annals.2005.162.1065}.

\bibitem{Ker39} R.~Kershner, The Number of Circles Covering a Set, {\it American J. of Math.}, {\bf 61}:3 (1939), 665--671, \url{https://archive.org/search?query=sim_pubid%3A1913+AND+volume%3A61&sort=date}.

\bibitem{KZ73} A. Korkine and G. Zolotareff, Sur les formes quadratiques, {\it Math. Annalen}, {\bf 6} (1873), 366--389, \url{http://eudml.org/doc/156622}.

\bibitem{KZ77} A. Korkine and G. Zolotareff, Sur les formes quadratiques positives, {\it Math. Annalen}, {\bf 11} (1877), 242--292, \url{https://eudml.org/doc/156750}.

\bibitem{Min97} H.~Minkowski, Allgemeine Lehrsätze über die convexen Polyeder, {\it Nachr. Ges. Wiss. Gött., Math.-Phys. Kl.}, (1897), 198--220.

\bibitem{Mus08} O.~R.~Musin, The kissing number in four dimensions, {\it Ann. of Math.}, {\bf 168}:1 (2008), 1--32, \url{https://doi.org/10.4007/annals.2008.168.1}.

\bibitem{PZ04} F.~Pfender, G.~M.~Ziegler, Kissing numbers, sphere packings, and some unexpected proofs,
{\it Notices AMS}, {\bf 51}:8 (2004), 873--883. \url{https://www.ams.org/journals/notices/200408/fea-pfender.pdf}.

\bibitem{Rog59} C.~A.~Rogers, Lattice coverings of space, {\it Mathematika}, {\bf 6}:1 (1959), 33--39, \url{https://doi.org/10.1112/S002557930000190X}.

\bibitem{RB78} S.~S.~Ryshkov and E.~P.~Baranovskii, $C$-Types of $n$-Dimensional Lattices and 5-Dimensional Primitive Parallelohedra (with Application to the Theory of Coverings), American Mathematical Society, 1978.

\bibitem{Sch09} A.~Sch\"urmann, Computational Geometry of Positive Definite Quadratic Forms, {\it Univ. Lecture Ser.}, {\bf 48}, AMS, 2009.

\bibitem{SV05} A.~Sch\"urmann and F.~Vallentin, Local covering optimality of lattices: Leech lattice versus root lattice $\mathsf E_8$, {\it IMRN}, 32 (2005), 1937--1955, \url{https://doi.org/10.1155/IMRN.2005.1937}.


\bibitem{SV06} A.~Sch\"urmann and F.~Vallentin, Computational Approaches to Lattice Packing and Covering Problems, {\it Disc. \& Comp. Geom.}, {\bf 35}:1 (2006), 73--116, \url{https://doi.org/10.1007/s00454-005-1202-2}.


\bibitem{Smi88} W.~D.~Smith, Ph.D. thesis: Studies in computational geometry motivated by mesh generation,
Department of Applied Mathematics, Princeton University, 1988.

\bibitem{Val03} F.~Vallentin, Sphere Coverings, Lattices, and Tilings (in Low Dimensions), PhD dissertation, TU Munich, 2003, \url{https://mediatum.ub.tum.de/doc/602017/document.pdf}.

\bibitem{Via17} M.~S.~Viazovska, The sphere packing problem in dimension 8, {\it Ann. of Math.}, {\bf 185}:3 (2017), 991--1015, \url{https://doi.org/10.4007/annals.2017.185.3.7}.

\bibitem{Vor08} G.~Voronoi, Nouvelles applications des param\`etres continus \`a la th\'eorie des formes quadratiques. {\it Crelle's Journal}, {\bf 133} (1908), 97--178, \url{https://doi.org/10.1515/crll.1908.133.97}; {\bf 134} (1908), 198--287, \url{https://doi.org/10.1515/crll.1908.134.198}; and {\bf 136} (1909), 67--181, \url{https://doi.org/10.1515/crll.1909.136.67}. 


\bibitem{Zong99} C.~Zong, Sphere packings, Springer-Verlag, 1999, \url{https://doi.org/10.1007/b98975}.
\end{thebibliography}
\end{document}